\documentclass[a4paper,11pt]{article}
\usepackage{amssymb}
\usepackage{bm}
\usepackage{amsmath}
\usepackage{graphicx}
\usepackage{xcolor}
\usepackage{dsfont}
\usepackage{caption}

\def\vsp{\vspace{0.2em}\noindent}
\def\vspp{\vspace{0.5em}\noindent}
\def\vsppp{\vspace{0.8em}\noindent}
\def\espace{\vspace{1.2em} }

\def\proof{\noindent {\bf Proof \ }}
\def\eproof{{ }\hfill$\Box$}

\def\noi{\noindent}
\def\longrightarrow{\to}

\def\a{\alpha}
\def\b{\beta}

\def\greta{\pmb\eta}

\def\RR{{\mathbb R} }
\def\CC{{\mathbb C}}
\def\NN{{\mathbb N}}
\def\ZZ{{\mathbb Z}}
\def\HH{{\mathbb H}}
\def\DD{{\mathbb D} }
\def\AA{{\mathbb A} }
\def\SS{{\mathbb S} }

\def\hh{{\cal H} }
\def\vv{{\cal V} }

\def\C{{\cal C} }

\def\H{{\cal H} }
\def\D{{\cal D}}

\def\ci{c_\circ} 
\def\pX{\partial_\infty X}
\def\conv{{\rm conv} }

\def\Im{{\rm Im\, }}
\def\is{{\rm Is\, }}
\def\plu{[0,\infty[}
\def\inj{{\rm inj}}
\def\join{{\rm join}\, }

\newcommand{\plouf}{\tag*{$\Box$}}

\newtheorem{theo}{Theorem} [section]
\newtheorem{prop}[theo]{Proposition}
\newtheorem{lemme}[theo]{Lemma}

\newtheorem{cor}[theo]{Corollary}

\newtheorem{defi}[theo]{Definition}
\newtheorem{notation}[theo]{Notation}
\newtheorem{fact}[theo]{Fact}

\newtheorem{hyp}[theo]{Hypothesis}

\newtheorem{example}[theo]{Example}
\newtheorem{remark}[theo]{Remark}
\numberwithin{equation}{section}

\def\bfa{\begin{fact}}
	\def\efa{\end{fact}}

\def\bhyp{\begin{hyp}}
	\def\ehyp{\end{hyp}}

\def\bt{\begin{theo}}
\def\et{\end{theo}}

\def\bc{\begin{cor}}
\def\ec{\end{cor}}

\def\bd{\begin{defi}}
\def\ed{\end{defi}}

\def\bl{\begin{lemme}}
\def\el{\end{lemme}}

\def\bp{\begin{prop}}
\def\ep{\end{prop}}

\def\bex{\begin{example}}
\def\eex{\end{example}}

\def\bre{\begin{remark}}
\def\ere{\end{remark}}

\def\bno{\begin{notation}\rm}
	\def\eno{\end{notation}}

\begin{document}
\title{Harmonic quasi-isometric  maps III~:\\  
quotients of Hadamard manifolds}
\author{Yves Benoist \&
Dominique Hulin}
\date{}

\vfill\eject
\setcounter{page}{1}
\maketitle

\begin{abstract}
In a previous paper, we proved that a quasi-isometric map $f:X\longrightarrow Y$ between two pinched Hadamard manifolds $X$ and $Y$ is within bounded distance from a unique harmonic map.

 We extend this result to maps $f:\Gamma\backslash X\longrightarrow Y$, where $\Gamma$ is a convex cocompact  discrete group of isometries of $X$ and $f$ is locally quasi-isometric at infinity.
\end{abstract}

\renewcommand{\thefootnote}{\fnsymbol{footnote}} 
\renewcommand{\thefootnote}{\arabic{footnote}} 

\section{Introduction}
\subsection{Statement and history}\label{subsec: statement}
The main result in this paper, which is a continuation of \cite{BH15}, \cite{BH17} and \cite{BH18}, is the following.
\bt\label{th main} Let $X$ and $Y$ be pinched Hadamard manifolds and $\Gamma\subset \is (X)$ be a torsion-free convex cocompact discrete subgroup of the group of isometries of $X$. Assume that the quotient manifold $M=\Gamma\backslash X$ is not compact.

\vspace{0.2em} Let $f:M\longrightarrow Y$  be a map. Assume that $f$ is quasi-isometric or, more generally, that $f$ is locally quasi-isometric at infinity (see Definition \ref{def loc qi}).

\vspace{0.2em} Then, there exists a unique harmonic map $h:M\longrightarrow Y$  within bounded distance from $f$, namely such that
$\displaystyle d(f,h):=\sup_{m\in M}d(f(m),h(m))<\infty$.
\et

\vspace{0.6em}
Let us begin with a short historical background (see  \cite[Section 1.2]{BH15} for more references). In the 60's, Eells and Sampson prove in \cite{EellsSampson64} that any smooth map $f:M\longrightarrow N$ between compact Riemannian manifolds, where $N$ is assumed to have non positive curvature, is homotopic to a harmonic map $h$. This harmonic map $h$  actually minimizes the Dirichlet energy $\int_M |Dh|^2 $ among all maps that are homotopic to $f$.

Later on,  P.Li and J.Wang conjecture in \cite{LiWang98} that it is possible to relax the co-compactness assumption in the Eells-Sampson theorem~: namely, they conjecture that any quasi-isometric map $f:X\longrightarrow Y$ between non compact rank one symmetric spaces is within bounded distance from a unique harmonic map. This  extends a former conjecture by R. Schoen in \cite{Schoen90}.
The Schoen-Li-Wang conjecture is proved by Markovic and Lemm-Markovic (\cite{Marko2}, \cite{Marko3}, \cite{Markon})
when  $X=Y=\HH^n$ are a real hyperbolic space, and in our papers 
\cite{BH15}, \cite{BH18} when $X$ and $Y$ are either rank one symmetric spaces or, more generally, pinched Hadamard manifolds. See also the recent paper by Sidler and Wenger \cite{SidWen18} and the survey by Guéritaud \cite{Gueritaud19}.
Our theorem \ref{th main} generalizes these results by allowing some topology in the source manifold.

\subsection{Examples and definitions}\label{subsec:examples and derfinitions}
A first concrete example where our theorem applies is the following, which is illustrated in Figure (a).
\bc\label{cor 1 } Let $\Gamma\subset {\rm PSl}_2\RR$ be a convex cocompact Fuchsian group. 

Any quasi-isometric map $f:\Gamma \backslash \HH^2\to \HH^3$ is within bounded distance from a unique harmonic map $h:\Gamma \backslash \HH^2\to \HH^3$.
\ec
\begin{figure}[h]
	\includegraphics[width=12.5cm]{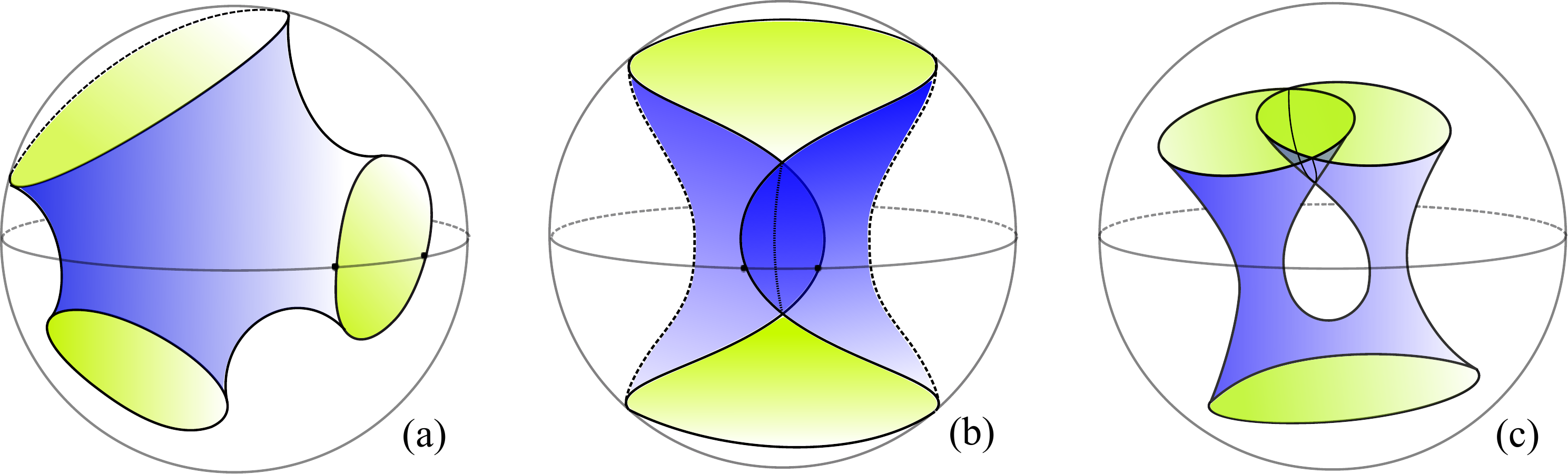}
	\hspace*{-3cm}\caption*{\small Examples of harmonic maps from surfaces to $\HH^3$,  with prescribed boundary values at infinity, such as  discussed in this paper.}
\end{figure}
As we will see later on  in Paragraph \ref{subsec: gromov hyp}, proving Theorem \ref{th main} amounts to solving a Dirichlet problem at infinity. Recall that a map $h_\infty :\SS^{k-1}\to \SS^{n-1}$ ($k,n\geq 2$) is quasi-regular if it is locally the boundary value of some quasi-isometric map from $\HH^k$ to $\HH^n$.
The following concrete example is again a special case of our theorem, illustrated in Figure (b).
\bc\label{cor 2 localement qi}
Let  $h_\infty :\SS^{k-1}\to \SS^{n-1}$ be a quasi-regular map. Then, $h_\infty$ extends as a  harmonic map 
$h:\HH^k\to \HH^n$. 
\ec
When $k=n$, Corollary \ref{cor 2 localement qi} was proved by  Pankka and Souto in \cite{PanSou17}.

Figure (c) illustrates Theorem \ref{th main} in a situation that combines those of Corollaries \ref{cor 1 } and \ref{cor 2 localement qi}.

\vspace{0.8em}
Let us now explain the hypotheses and conclusion of Theorem \ref{th main}.
A map $f:X\longrightarrow Y$ between two metric spaces is quasi-isometric if there exists a constant $c\geq 1$ such that 
\begin{align}\label{def qi}
c^{-1\, }d(x,x')-c\leq d(f(x),f(x'))\leq c\, d(x,x')+c
\end{align}
holds for any $x,x'\in X$. Such a map needs not be continuous.

A smooth map $h:X\longrightarrow Y$ between Riemannian manifolds is harmonic if it is a critical point for the Dirichlet energy $\int |Dh|^2$, namely if it satisfies the elliptic PDE 
$$
{\rm Tr}\, D^2h=0\, .
$$ 

A pinched Hadamard manifold is a complete simply-connected Riemannian manifold $X$ with dimension at least 2 whose sectional curvature is pinched between two negative constants, namely
\begin{equation}\label{eq pinching}
-b^2\leq K_X\leq -a^2<0\, .
\end{equation}
Observe that working with pinched Hadamard manifolds provides a natural and elegant framework to deal simultaneously with all the rank one symmetric spaces.

\vspace{0.5em}
A discrete subgroup $\Gamma\subset\is (X)$ is convex cocompact when the convex core $K\subset M$ is compact, see Definition \ref{def K}.
We will see  in Section \ref{sec:cvx cpct} (Proposition \ref{prop cvx ccpct}) that requiring  the discrete subgroup $\Gamma\subset\is (X)$ to be convex cocompact is equivalent to assuming that the Riemannian quotient $M=\Gamma\backslash X$ is Gromov hyperbolic, and that its injectivity radius ${\rm inj} :M\to\plu$ is a proper function on $M$. Therefore we can speak of the boundary at infinity $\partial_\infty M$ of $M$.
\bd\label{def loc qi} A map $f:M\to Y$ is locally quasi-isometric at infinity if it
satisfies the two conditions~:

{\rm (a) } the map $f$  is rough Lipschitz~: there exists a constant $c\geq 1$ such that  $d(f(m),f(p))\leq c\, d(m,p)+c$ for any $p,m\in M$~;

{\rm (b) }  each point  $\xi\in\partial_\infty M$ in the boundary at infinity of $M$ admits a neighbourhood $U_\xi$ in $M\cup\partial_\infty M$ such that the restriction $f:U_\xi\cap M\to Y$ is quasi-isometric.
\ed
Note that condition (a) follows from condition (b) when the map $f$ is  bounded on compact subsets of $M$.

In Section 3, we give a few counter-examples that emphasise the relevance of assuming in Theorem \ref{th main} that $\Gamma$ is convex cocompact.

Note that the assumption that $\Gamma$ is torsion-free in Theorem \ref{th main} is only used to ensure that the quotient $M=\Gamma\backslash X$ is a manifold. See Kapovich \cite{Kapovich18} for nice examples of non torsion-free groups.

\subsection{Structure of the proof}\label{subsec: structure}
Although the rough structure of the proof of Theorem \ref{th main} is similar to those of \cite{BH15} or \cite{BH18}, we now have to deal with new  issues. Two new crucial steps will be  understanding the geometry of the source manifold $M=\Gamma\backslash X$ in Proposition \ref{prop cvx ccpct}, and obtaining uniform estimates for the harmonic measures on $M$ in Corollary \ref{cor upper bound harmonic}.

Let us now give an overview of the proof of existence in Theorem \ref{th main}. We refer to Section \ref{sec:proof} for complete proofs of existence and uniqueness. 

Replacing the quasi-isometric map $f$ by local averages in Lemma \ref{lem smooth f}, we may assume that $f$ is smooth. This averaging process even allows us to assume that $f$ has uniformly bounded first and second covariant derivatives, a fact that will be crucial later on to prove the so-called boundary estimates of Paragraph \ref{subsec:boundary estimates}.

To prove existence, we begin by solving a family of  Dirichlet problems on bounded domains of $M $. Namely, we introduce in Proposition \ref{prop smooth core} an exhausting and increasing family of compact convex domains  $V_R\subset M $ with smooth boundaries ($R>0$), and consider for each $R>0$ the unique harmonic map $h_R:V_R\longrightarrow Y$ which is solution of the Dirichlet problem ``$h_R=f$ on the boundary $\partial V_R$'' (Lemma \ref{lem hR}).

The heart of the proof, in  Proposition \ref{prop rho}, consists in showing that the distances $d(h_R,f)$ are uniformly bounded 
by a constant $\bar\rho$ that does not depend on $R$.
Once we have this bound, we recall in Paragraph \ref{subsec:existence} a standard compactness argument that ensures that the family of harmonic maps $(h_R)$ converges, when $R$ goes to infinity,  to a harmonic map $h:M\longrightarrow Y$. The limit harmonic map $h$   still satisfies $d(h,f)\leq\bar\rho$.

\vspace{0.5em} 
Let us now briefly  explain how we obtain this uniform bound for the distances $\rho_R:=d(h_R,f)$. We assume by contradiction that $\rho_R$ is very large. 

The first step  consists in proving  that   the distance $\rho_R$ is achieved at a point $m\in M $ which is far away from the boundary $\partial V_R$. This  relies on the uniform bound for the covariant derivatives $Df$ and $D^2f$, and on the equality $d(h_R(p),f(p))=0$ when $p$ lies in the boundary $\partial V_R$. 

The second step consists in the so-called interior estimates of Section \ref{sec:interior estimates}. Since the point $m\in V_R$ such that $\rho_R=d(h_R(m),f(m))$ is far away from the boundary $\partial V_R$, we may select a large neighbourhood $W(m)\subset M $ of the point $m$, such that $W(m)\subset V_R$. 

Let us describe the neighbourhoods $W(m)$. We fix  two large constants $1\ll \ell_0 \ll\ell$  that will not depend on $R$ and need to be properly chosen. If $m$ is far enough from the convex core  $K\subset M $, the injectivity radius at the point $m$ will be larger than $\ell_0$ and  we will choose $W(m)$ to be the Riemannian ball $W(m)=B(m,\ell_0)\subset V_R$ with center $m$ and radius $\ell_0$. If $m$ is close to the convex core, we will chose $W(m)=V_\ell\subset V_R$ to be a fixed large compact neighbourhood of the convex core $K$. 

\vspace{0.2em}
To wrap up the proof in Section \ref{sec:interior estimates}, we  will exploit the subharmonicity of the function $q\in W(m)\longrightarrow d(f(m),h_R(q))\in\plu$. The crucial tool for obtaining a contradiction, and thus proving that $\rho_R$ cannot be very large, is uniform estimates for the harmonic measures of the boundaries of these  neighbourhoods $W(m)$ relative to the point $m$ that are proved in Sections \ref{sec:harmonic measure, lower} and \ref{sec:upper bound}. See Proposition \ref{prop lower bound harmonic} and Corollary \ref{cor upper bound harmonic}. We refer to Paragraph \ref{subsec: overview close} for a more precise overview of this part the proof.

\section{Convex cocompact subgroups of $\is (X)$}\label{sec:cvx cpct}

\begin{quote}
	In this section  we characterize convex cocompact subgroups  $\Gamma\subset \is (X)$ in terms of the geometric properties of the quotient $\Gamma\backslash X$.
\end{quote}

\subsection{Gromov hyperbolic metric spaces}\label{subsec: gromov hyp}
We first recall a few facts and definitions concerning Gromov hyperbolic metric spaces, see \cite{GhysHarp90}.

Let $\delta >0$. A geodesic metric space $X$ is said to be $\delta$-Gromov hyperbolic when every geodesic triangle $\Delta$ in $X$ is $\delta$-thin, namely when each edge of $\Delta$ lies in the $\delta$-neighbourhood of the union of the other two edges. For example, any pinched Hadamard manifold is Gromov hyperbolic for some constant $\delta_X$ that depends only on the upper bound for the curvature.

When the Gromov hyperbolic space $X$ is proper, namely its  balls are compact, the boundary at infinity $\partial_\infty X$ of $X$ may be defined as the set of equivalence classes of geodesic rays, where two geodesic rays are identified when they remain within bounded distance from each other. In case $X$ is a pinched Hadamard manifold, the boundary at infinity (or visual boundary) $\partial_\infty X$ naturally identifies with the tangent sphere  at any point $x\in X$. The boundary at infinity provides a compactification $\overline X =X\cup\partial_\infty X$ of $X$.

A quasi-isometric map $f:X\longrightarrow Y$ between proper Gromov hyperbolic spaces admits a boundary value $\partial_\infty  f:\partial_\infty X\to\partial_\infty Y$, and two quasi-isometric maps $f_1,f_2 : X\longrightarrow Y$ have the same boundary value if and only if $d(f_1,f_2)<\infty$.

\vspace{0.6em}
Since our quotient $M=\Gamma\backslash X$ is Gromov hyperbolic (see Proposition \ref{prop cvx ccpct}), we may thus rephrase the conclusion of Theorem \ref{th main} in case $f$ is assumed to be globally quasi-isometric.
\bc 	Given a quasi-isometric map $f:M\to Y$, there exists a unique harmonic   quasi-isometric map $h:M\longrightarrow Y$ which is a solution to the Dirichlet problem at infinity with boundary value 
\begin{equation*}\label{eq boundary pb at infinity}
\partial_\infty  h=\partial_\infty  f:\partial_\infty M\longrightarrow\partial_\infty Y\, .
\end{equation*} 
\ec

\subsection{Gromov product}\label{subsec: Gromov}
In Section \ref{sec:interior estimates}, we shall use Gromov products in the Gromov hyperbolic manifolds $ M$ and  $Y$ to carry out the proof of Theorem \ref{th main}. Here, we recall the definition and two properties of the Gromov product.

\vspace{0.2em}
In a metric space, the Gromov product of the three points $x,x_1,x_2$  is defined as
\begin{equation*}
(x_1,x_2)_x=\frac 12 (d(x,x_1)+d(x,x_2)-d(x_1,x_2))\, .
\end{equation*}
Gromov hyperbolicity may be expressed in terms of Gromov products. Also, in a Gromov hyperbolic space, the Gromov product $(x_1,x_2)_x$ is roughly equal to the distance from $x$ to a minimizing geodesic segment $[x_1,x_2]$.  
In particular, the following holds.
\bl\label{lem proprietes Gromov prdct}{\rm\cite[Chap.2]{GhysHarp90}} 
Let $X$ be a $\delta$-Gromov hyperbolic space. Then, for any points $x,x_1,x_2,x_3\in X$ and any   minimizing geodesic segment $[x_1,x_2]\subset X$, one has
\begin{align}   
(x_1,x_3)_x&\geq \min ((x_1,x_2)_x,(x_2,x_3)_x)-2\delta \label{eq gromov product inegalite triangulaire} \\
(x_1,x_2)_x&\leq d(x, [x_1,x_2])\leq (x_1,x_2)_x+2\delta\ \label{eq gromov product distance segment}.
\end{align}
\el
Moreover,  the next lemma tells us that Gromov products are quasi-invariant under  quasi-isometric maps.

\bl\label{lem gromov product quasi iso}{\rm \cite[Prop.5.15]{GhysHarp90}} Let $X$ and $Y$ be Gromov hyperbolic spaces, and $f:X\to Y$ be a $c$ quasi-isometric map. Then, there exists a constant $A>0$ such that,
for any three points $x,x_1,x_2\in X$~:
$$
c^{-1}\, (x_1,x_2)_x-A\leq (f(x_1),f(x_2))_{f(x)}\leq c\, (x_1,x_2)_x+A\, .
$$
\el

\subsection{Convex cocompact subgroups}

We begin with  definitions that are classical in the hyperbolic space $\HH^k$. See Bowditch  \cite{Bowditch95} for a reference dealing with pinched Hadamard manifolds.
\bd Let $X$ be a pinched Hadamard manifold and $\Gamma\subset\is (X)$ be a discrete subgroup of the group of isometries of $X$.

The limit set  $\Lambda_\Gamma\subset \partial_\infty X$ of the group $\Gamma$ is the closed subset of $\partial_\infty X$ defined as
$\Lambda_\Gamma = \overline{\Gamma x}\smallsetminus \Gamma x$, where $x$ is any point in $X$ and the closure $\overline{\Gamma x}$ of its orbit is taken in the compactification $\overline X = X\cup\partial_\infty X$.

The domain of discontinuity of $\Gamma$ is $\Omega_\Gamma = \partial_\infty X\smallsetminus \Lambda_\Gamma$. It is an open subset of the boundary at infinity $\partial_\infty X$. The group $\Gamma$ acts properly discontinuously on $X\cup\Omega_\Gamma\subset\overline X$. 
\ed
Recall that a subset $C\subset M$ of a Riemannian manifold is geodesically convex if, given two points  $m_1,m_2\in C$, any minimizing geodesic segment $[m_1,m_2]$ joining these two points  lies in $C$.

\bd\label{def K} Let $X$ be a pinched Hadamard manifold, $\Gamma\subset\is (X)$ be an infinite discrete subgroup of the isometry group of $X$ and $M=\Gamma\backslash X$. 

The convex hull $\conv\Lambda_\Gamma\subset X$ of $\Lambda_\Gamma$ is the smallest closed geodesically convex subset $C$ of $X$ such that $\overline C\smallsetminus C=\Lambda_\Gamma$ (where  the closure $\overline C$ is again taken in the compactification $\overline X$).
The convex core of $M$ is the quotient $K:=\Gamma \backslash (\conv\Lambda_\Gamma)\subset M$.

The group $\Gamma$ is said to be convex cocompact if the convex core $K\subset M$ is compact.
This is equivalent to requiring the quotient $\Gamma\backslash (X\cup \Omega_\Gamma)\subset \Gamma\backslash \overline X$ to be compact.  See \cite[Th.6.1]{Bowditch95}.
\ed

The hypothesis in  Theorem \ref{th main} that $\Gamma$ is convex cocompact will be used in this paper through the following characterization in terms of geometric properties of the Riemannian quotient $M$.
\bp\label{prop cvx ccpct}
Let $X$ be a pinched Hadamard manifold and $\Gamma\subset \is (X)$ be an infinite torsion-free discrete subgroup. Then, the group $\Gamma$ is convex cocompact if and only if $M=\Gamma\backslash X$ satisfies the following conditions~:
\begin{itemize}
	\item the quotient $M$ is Gromov hyperbolic~;
	\item the injectivity radius ${\rm inj}: M \to \plu$ is a proper function.
\end{itemize}
\ep
Proposition \ref{prop cvx ccpct} is proved in Paragraphs \ref{sub : G hyp quotient}, \ref{sub : inj rad quotient} and \ref{sub : converse} below.

\subsection{Gromov hyperbolicity of the quotient $M=\Gamma\backslash X$} \label{sub : G hyp quotient}
In this paragraph, $X$ is our pinched Hadamard manifold, $\Gamma\subset\is (X)$ is a torsion-free convex cocompact subgroup and $M=\Gamma\backslash X$.
We want to prove that there exists a constant $\delta >0$ such that every triangle $T\subset M$ in the quotient is $\delta$-thin. To do this, we start with quadrilaterals. 
We first recall a result by Reshetnyak that compares quadrilaterals in a Hadamard manifold with quadrilaterals in the model hyperbolic plane $\HH^2(-a^2)$ with constant curvature $-a^2$. 
\bl{\small\bf Reshetnyak  comparison lemma} {\rm \cite{Reshetnyak68}}\label{lem resh}

\vsp
Let $X$ be a Hadamard manifold satisfying the pinching assumption \eqref{eq pinching}.
For every quadrilateral $[x,y,z,t]$ in $X$, there exists a convex quadrilateral $[\bar x, \bar y,\bar z,\bar t]$ in  $\HH^2(-a^2)$ and a  map
$j:[\bar x, \bar y,\bar z,\bar t]\to X$ that is $1$-Lipschitz, that sends respectively the vertices $\bar x, \bar y,\bar z,\bar t$ on $x,y,z,t$, and whose restriction to each one of the four edges of $[\bar x, \bar y,\bar z,\bar t]$  is isometric.
\el
We now deduce from the Reshetnyak Lemma a standard property of quadrilaterals in $X$, that will be used again in the proof of Proposition \ref{prop tau0}. The
Gromov hyperbolicity of $X$ ensures that any edge of a quadrilateral in $X$ lies in the $2\delta_X$-neighbourhood of the union of the three other edges. One can say more under an angle condition. Namely~:

\bl{\small\bf Thin quadrilaterals in $X$}\label{lem quad X} \\
Let $X$ be a pinched Hadamard manifold. 
Let $\varepsilon >0$. There exists an angle $0<\alpha_\varepsilon <\pi/2$ and a distance $\tilde\delta_\varepsilon $ such that for any quadrilateral $[x,x_1,y_1,y]$ in the Hadamard manifold $X$ with $d(x_1,y_1)\geq\varepsilon$, and whose angles at both vertices $x_1$ and $y_1$ satisfy $\angle_{x_1}\geq\pi/2-\a_\varepsilon $ and $\angle_{y_1}\geq\pi/2-\a_\varepsilon $, the $\tilde\delta_\varepsilon$-neighbourhood of the edge $[x,y]$ contains the union of the other three edges, namely~: 
$$
[x,x_1]\cup[x_1,y_1]\cup[y_1,y]\subset \vv_{\tilde\delta_\varepsilon}([x,y])\, .
$$
\el
\proof When $X$ is the hyperbolic plane, the easy proof is  left to the reader. 
The general case follows  from this special case and the Reshetnyak comparison lemma \ref{lem resh}. Indeed,  the comparison quadrilateral   also satisfies the distance and angles conditions at the points $\bar x_1$ and $\bar y_1$.
\eproof

\espace We now turn to quadrilaterals in our quotient space $M =\Gamma\backslash X$. In the whole paper, unless otherwise  specified, all triangles and quadrilaterals in $M$ will be assumed to be minimizing, namely their edges will be minimizing geodesic segments.
\bl{\small\bf Thin quadrilaterals in $M $}\label{lem quad M}

\vsp
Assume that $X$ is a pinched Hadamard manifold, that $\Gamma\subset \is (X)$ is a torsion-free convex cocompact subgroup, and let $M=\Gamma\backslash X$.
Let $V\subset M $ be a compact convex subset of $M$ whose lift $\tilde V\subset X$ is convex. 

(1) There exists $\delta_0$ such that, for any  quadrilateral $[p,p_1,q_1,q]$ in $M $ with both $p_1,q_1\in V$,  any edge of this quadrilateral lies in the $\delta_0$-neighbourhood of the union of the three other edges.

(2) Let $\varepsilon >0$. There exists  $\delta_\varepsilon$ such that if we assume moreover that $p_1,q_1\in V$ are respectively the projections of the points $p$ and $q$ on the convex set $V$, and that $d(p_1,q_1)\geq\varepsilon$, then~:
$$
[p,p_1]\cup[p_1,q_1]\cup[q_1,q]\subset \vv_{\delta_\varepsilon}([p,q])\, .
$$
\el
\proof Lift successively the minimizing geodesic segments  $(p_1p)$, $(pq)$ and $(qq_1)$ to geodesic segments  $(x_1x)$, $(xy)$ and $(yy_1)$  in $X$, so that the geodesic segment $(x_1y_1)$ projects to a curve $c$ from $p_1$ to $q_1$ that lies in $V$. 

Since both the curve $c$ and the geodesic segment lie in the compact set $V$, the conclusion follows from Lemma \ref{lem quad X}, with $\delta_0=2\delta_X+d_V$ and $\delta_\varepsilon=\tilde\delta_\varepsilon+d_V$, where $d_V$ denotes the diameter of $V$.
\eproof

\bc {\small\bf Gromov hyperbolicity of the quotient $\Gamma\backslash X$}

\vsp
Assume that $X$ is a pinched Hadamard manifold and that $\Gamma\subset \is (X)$ is a torsion-free convex cocompact subgroup. Then, the quotient $M=\Gamma\backslash X$ is Gromov hyperbolic.
\ec
\proof
Let $V\subset M $ be a compact convex neighbourhood of the convex core with smooth boundary, whose lift in $X$ is convex. Such a neighbourhood will be constructed in Proposition \ref{prop smooth core}. Let $\varepsilon={\rm inj} (V)/2$ (where ${\rm inj} (V)=\inf_{m\in V}\inj(m)$ denotes the injectivity radius on $V$).

\vspace{0.8em}
Let $T=[p,q,r]$ be a triangle in $M $.  
In   case where at least one of the vertices of $T$ lies in $V$,  Lemma \ref{lem quad M} applied to  a quadrilateral with two equal vertices proves that $T$ is $\delta_0$-thin.

\vspace{0.5em}
We now turn to the case where none of $p,q,r$ lie in $V$ and denote by $p_1,q_1,r_1$ their projections on $V$. 

Assume first that $d(p_1,q_1)\leq\varepsilon$ and $d(q_1,r_1)\leq \varepsilon$. Then, the triangle $[p_1,q_1,r_1]$ is homotopically trivial. Writing $p={\rm exp}_{p_1}u$ with $u$ a normal vector to $V$ at point $p_1$ (and similarly for $q$ and $r$) we construct an homotopy between $[p,q,r]$ and a constant map, so that the triangle $[p,q,r]$ lifts to a triangle $[x,y,z]$ in $X$. Since $X$ is $\delta_X$-Gromov hyperbolic, $[x,y,z]$ is $\delta_X$-thin, hence $[p,q,r]$ is $\delta_X$-thin too.

Assume now that $d(p_1,q_1)\geq\varepsilon$ and $d(q_1,r_1)\geq \varepsilon$. The first part of Lemma \ref{lem quad M}, applied to the quadrilateral $[p,p_1,r_1,r]$, yields that $[p,r]$ lies in the $\delta_0$-neighbourhood of $[p,p_1]\cup[p_1,r_1]\cup[r_1,r]$, so that  $[p,r]$ lies in the $(\delta_0+d_V)$-neighbourhood of $[p,p_1]\cup[r_1,r]$  (recall that $d_V$ is the diameter of $V$). The second part of Lemma \ref{lem quad M} now applied to both quadrilaterals $[p,p_1,q_1,q]$ and $[q,q_1,r_1,r]$ yields that 
$[p,r]$ lies in the $(\delta_0+\delta_\varepsilon+d_V)$-neighbourhood of $[p,q]\cup[q,r]$.

Assume finally that $d(p_1,q_1)\leq\varepsilon$ and $d(q_1,r_1)\geq \varepsilon$. It follows from the first part of Lemma \ref{lem quad M}  that $[p,r]$ lies in the $(\delta_0+d_V)$-neighbourhood of $[p,p_1]\cup[r_1,r]$ and that $[p,p_1]$ lies in the $(\delta_0+d_V)$-neighbourhood of $[p,q]\cup[q_1,q]$, while the second part of this Lemma ensures that
$[q,q_1]\cup[r_1,r]$ lies in the $\delta_\varepsilon$-neighbourhood of $[q,r]$. 

Hence, every triangle in $M$ is $(2(\delta_0+d_V)+\delta_\varepsilon+\delta_X)$-thin, so that $M$ is Gromov hyperbolic.
\eproof

\subsection{Injectivity radius of the quotient $M=\Gamma\backslash X$}\label{sub : inj rad quotient}
Our aim in this paragraph is the following.
\bp{\bf\small Properness of the injectivity radius}\label{prop propre inj radius}

\vsp Assume that $X$ is a pinched Hadamard manifold and that $\Gamma\subset \is (X)$ is a torsion-free convex cocompact subgroup. Let $M=\Gamma\backslash X$.

Then, the injectivity radius  ${\rm inj} : M\to\plu$ is a proper function.
\ep

\proof For $p\in M $, we denote by $i(p)$ the injectivity radius at the point $p$, and let $j(p)=\inf\{ d(\tilde p,\gamma \tilde p) \, |\, \gamma\in\Gamma^*\}$, where $\tilde p\in X$ is a lift of $p$ and $\Gamma^*$ is the set of non trivial elements in $\Gamma$. Since $M $ has non positive curvature, there are no conjugate points in $M $ hence $j(p)=2i(p)$.

We proceed by contradiction and assume that there exists a sequence $(p_n)$ of points in $M $ going to infinity, and such that the injectivity radii $i(p_n)$ remain bounded. Let $q_n$ denote the projection of the point $p_n$ on the convex core $K\subset M $ ($n\in\NN$). Since $\Gamma$ is convex cocompact, there exists a compact $L\subset X$ such that each $q_n$ lifts in $X$ to a point $\tilde q_n\in L$.
Then, the geodesic segment $[q_n,p_n]$ lifts as $[\tilde q_n,\tilde p_n]$. By hypothesis, there exists $I>0$ and a sequence $\gamma_n\in\Gamma^*$
such that $d(\tilde p_n,\gamma_n\tilde p_n)\leq I$. Since the projection $X\to{\rm conv}(\Lambda_\Gamma)$ commutes to the action of $\Gamma$ and is $1$-Lipschitz, it follows that $d(\tilde q_n,\gamma_n\tilde q_n)\leq I$. By compactness of $L$, and since $\Gamma $ is discrete, one may thus assume that the sequence $(\tilde q_n)$ converges to a point $\tilde q\in L$ and that the bounded sequence $(\gamma_n)$ is constant, equal to $\gamma\in \Gamma^*$. The boundary at infinity $\partial_\infty X$ being compact, we may also assume that the sequence of geodesic rays $([\tilde q_n,\tilde p_n\mathclose[)$ converges to a geodesic ray $[\tilde q,\xi\mathclose[$ where $\xi\in \partial_\infty X$. 

By construction, the geodesic ray $[\tilde q,\xi\mathclose[$ is within bounded distance $I$ from its image $[\gamma\tilde q,\gamma\xi\mathclose[$, hence $\xi\in \partial_\infty X$ is a fixed point of $\gamma$. The group $\Gamma$ being discrete and torsion-free, it has no elliptic elements, so that $\xi\in\Lambda_\Gamma$. Thus, the whole geodesic ray $[\tilde q,\xi\mathclose[$ lies in ${\rm conv} (\Lambda_\Gamma)$, a contradiction to the fact that the sequence $(p_n)$ goes to infinity in $M $.
\eproof

\subsection{Converse}\label{sub : converse}

In this paragraph, we complete the proof of Proposition \ref{prop cvx ccpct} by proving the following.
\bp\label{prop converse}
Let $X$ be a pinched Hadamard manifold, and $\Gamma\subset \is (X)$ be a torsion-free discrete subgroup.
Assume that the quotient manifold $M=\Gamma\backslash X$ is Gromov hyperbolic, and that the injectivity radius $\inj :M\to\plu$ is a proper function.

Then, the group $\Gamma$ is convex cocompact.
\ep
We want to prove that the convex core $K:=\Gamma \backslash (\conv\Lambda_\Gamma)$ is a compact subset of $M$, where $\conv\Lambda_\Gamma$ denotes the convex hull of the limit set of $\Gamma$. We will rather work with the join of the radial limit set $\Lambda_\Gamma^r$.
\bd 
A geodesic ray $c:\plu \to M$ is said to be recurrent when it is not a proper map, that is if there exists a sequence $t_n\to +\infty$ and a compact set $L_c\subset M$ (that might depend on $c$) with $c (t_n)\in L_c$.

The radial limit set $\Lambda_\Gamma^r\subset\Lambda_\Gamma$ of $\Gamma$ is the set of endpoints $\xi\in \partial_\infty X$ of  geodesic rays in $X$ that project to a recurrent geodesic ray in $M$.
\ed
Since there exist closed geodesics in the negatively curved manifold $M$, the radial limit set $\Lambda_\Gamma^r$ is not empty, and it is a $\Gamma$-invariant subset of $\pX$. Hence, the limit set $\Lambda_\Gamma$ being a minimal $\Gamma$-invariant closed subset of $\partial_\infty X$, it follows that $\Lambda_\Gamma^r$ is dense in $\Lambda_\Gamma$.

\bd The join of a closed subset $Q\subset \partial_\infty X$ is defined as
$$
\join Q=\bigcup_{\xi_1,\xi_2\in Q}\mathopen]\xi_1,\xi_2\mathclose[\subset X\, 
$$
where $\mathopen]\xi_1,\xi_2\mathclose[$ denotes the geodesic line with endpoints $\xi_1$ and $\xi_2$.
\ed
The join is thus the first step towards the construction of the convex hull. One has $\join Q\subset\conv Q$. The following result by Bowditch \cite{Bowditch95} investigates how much we miss in the convex hull by considering only the join.
\bp {\rm\cite[Lemma 2.2.1 and Proposition 2.5.4]{Bowditch95}}\label{prop bowditch}

\vsp Let $X$ be a pinched Hadamard manifold.
There exists a real number $\lambda>0$ that depends only on the pinching constants such that, for any closed subset $Q\subset\partial_\infty X$, the convex hull of $Q$ lies in the $\lambda$-neighbourhood of the join of $Q$.
\ep
Let us now proceed with our proof. 
The  projection $K^r_j:=\Gamma\backslash(\join \Lambda_\Gamma^r)\subset M$ of the join of the radial limit set is the union of all the recurrent geodesics on $M$, namely of all geodesics that are recurrent both in the future and in the past.  We  begin with the following.
\bp\label{prop converse join} Under the assumptions of Proposition \ref{prop converse}, the join of the radial limit set of $\Gamma$  projects in $M$ to a bounded set $K^r_j$.
\ep
\proof Proposition \ref{prop converse join} is an immediate consequence of the definition of the join and of Lemma \ref{lem join borne} below.
\eproof

\vspace{0.6em}\noi
For  $r>0$, define $M_r=\{ m\in M \, |\; \inj (m)\leq r\}$. Under our hypotheses, $M_r$ is a compact subset of $ M$. Let $\delta>0$ such that  $ M$ is $\delta$-Gromov hyperbolic.

\bl\label{lem minimisante} Any geodesic segment in $  M$   whose image lies outside the compact subset $M_{3\delta}\subset M$ is minimizing.
\el
\proof We proceed by contradiction and assume that the geodesic segment $c:[0,l]\to M\smallsetminus M_{3\delta}$ is minimizing, but ceases to be minimizing on any larger interval. Since $ M$ has negative curvature, there are no conjugate points, so that there exists another minimizing geodesic $ \bar c:[0,l]\to  M$ with the same endpoints $c_0$ and $c_l$ as $c$.

The manifold $ M$ being $\delta$-hyperbolic, the geodesic segments $c$ and $\bar c$ are within distance $\delta$ from each other.
We may thus choose a subdivision $(x_p)_{0\leq p\leq P}$ of $c([0,l])$ with $x_0=c_0$, $x_P=c_l$, and such that $d(x_p,x_{p+1})\leq\delta$ when $0\leq p <P$, and another sequence $(y_p)_{0\leq p\leq P}$ with $y_p\in \bar c([0,l])$ and $d(x_p,y_p)\leq\delta$  when $0\leq p \leq P$. In particular, $d(y_p,y_{p+1})\leq 3\delta$, so that the length of any of the quadrilaterals $q_p=[x_p,x_{p+1},y_{p+1},y_p]$ is at most $6\delta$. Since the injectivity radius at the point $x_p$ is larger than $3\delta$, it follows that each quadrilateral $q_p$ is homotopically trivial, so that $c([0,l])$ and $\bar c([0,l])$ themselves are homotopic. But $X$ being a Hadamard manifold then yields $c=\bar c$, which is a contradiction. \eproof

\bl\label{lem join borne} Let $d_{3\delta}$ denote the diameter of $M_{3\delta}$. 
Under the assumption of Proposition \ref{prop converse}, any recurrent geodesic in $ M$ lies in a fixed compact subset of $ M$. More precisely, such a geodesic lies in the $d_{3\delta}$-neighbourhood of $M_{3\delta}$. 
\el
\proof Let $c:\RR\to M$ be a geodesic that lifts to a geodesic $\tilde c:\RR\to X$ with both endpoints in $\Lambda^r_\Gamma$.

We first claim that both geodesic rays $c_{|[0,\infty[}$ and $c_{|]-\infty,0]}$ must keep visiting $M_{3\delta}$. If this were not the case, Lemma \ref{lem minimisante} would ensure that one of these geodesic rays would be eventually minimizing, thus contradicting the assumptions that both endpoints of $c$ lie in the radial limit set.

Now assume by contradiction that the image of $c$ does not lie in the $d_{3\delta}$-neighbourhood of $M_{3\delta}$. Then we can find an interval $[a,b]$ such that $c(a),c(b)\in M_{3\delta}$ but with $c(t)\notin M_{3\delta}$ for every $t\in ]a,b[$, and such that $c([a,b])$ exits the $d_{3\delta}$-neighbourhood of $M_{3\delta}$. Hence  this geodesic segment $c([a,b])$ has length at least $2d_{3\delta}$. By Lemma \ref{lem minimisante} it would be minimizing, a contradiction to the fact that $d(c(a),c(b))\leq d_{3\delta}$.
\eproof

\vspace{0.6em}\noi{\bf Proof of Proposition \ref{prop converse}} We  noticed earlier that the radial limit set $\Lambda^r_\Gamma\subset\Lambda_\Gamma$ is dense in the limit set of $\Gamma$. Thus, the join of the full limit set lies in the closure of the join of the radial limit set. Hence Proposition \ref{prop converse join} ensures that $\join \Lambda_\Gamma$ projects to a bounded subset of $ M$. 

Now Proposition \ref{prop bowditch} ensures that $\conv\Lambda_\Gamma$ also projects to a bounded subset of $ M$ or, in other words, that the group $\Gamma$ is convex cocompact. \eproof

\section{Examples}\label{sec:examples}
\begin{quote}
	Before going into the proof of Theorem  \ref{th main}, we proceed with a few examples and counter-examples.
\end{quote}
\subsection{Two examples on the annulus}\label{subsec: example}
We begin with a family of  straightforward applications of Theorem \ref{th main}.

\medskip
Let $\AA (r)=\{ z\in\CC \, |\; 1/r<|z|<r\}$ and $\DD=\{ z\in\CC \, |\; |z|<1\}$  be an annulus ($r>1$) or the disk  equipped with their complete hyperbolic metrics. 
\bex  For every $\alpha\in\RR$, there exists a unique harmonic quasi-isometric map $h_\a:\AA (r)\to\DD$ with boundary value  $g_\a :\partial_\infty (\AA (r))\to\partial_\infty \DD$ defined as
\begin{equation*}
g_\a(re^{i\theta})=e^{i\theta} \qquad  g_\a(e^{i\theta}/r)=e^{i(\theta +\a)}\, .
\end{equation*}
\eex
Note that the case where the parameter $\a $ is equal to $0$ is rather easy. Indeed, the uniqueness of the harmonic map with boundary value $g_0$ ensures that $h_0$ has a symmetry of revolution and thus reads as $h_0(te^{i\theta})=u(t)e^{i\theta}$, and the condition that $h_0$ is harmonic reduces to a second order ordinary differential equation  on $u$ whose solutions can be expressed in terms of elliptic integrals.

\espace\noi 
We now give an illustration of Theorem \ref{th main} where the restriction of the boundary map to each connected component of the boundary at infinity $\partial_\infty M$  is not injective.
\bex
The map  $g :\partial_\infty (\AA (r))\to\partial_\infty \DD$ defined as
\begin{equation*}
g(re^{i\theta})=e^{2i\theta} \qquad  g_\a(e^{i\theta}/r)=e^{-3i\theta }
\end{equation*}
is the boundary value of a harmonic map $h:\AA (r)\to\DD$.
\eex

\subsection{First counter-example~: surfaces admitting a cusp}\label{subsec:cusp}

In the next two paragraphs, we provide counter-examples to emphasize the importance of assuming that $\Gamma$ is convex cocompact. 

We first prove non-existence for hyperbolic surfaces with a cusp.
\bp\label{propwith cusp}
Let $S$ be a non compact hyperbolic surface of finite topological type admitting at least one cusp, and $Y$ be a pinched Hadamard manifold. Then, there exists no harmonic quasi-isometric map $h:S\to Y$.
\ep
Note that, such a surface $S$ being  quasi-isometric to the wedge of a finite number of hyperbolic disks and rays,   there always exist quasi-isometric maps $f:S\to \HH^n$ for any $n\geq 3$ -- none of them harmonic.  This proposition is an immediate consequence of the following lemma.
\bl\label{lem cusp}
Let $\tau >0$ and $  \Sigma $ be the quotient of  $\{ z\in\CC \, |\; {\rm Im}z\geq 1\}$,  equipped with the hyperbolic metric $\frac{|dz|^2}{(\Im z)^2}$ under the map $z\to z+\tau$. 

Let $Y$ be a pinched Hadamard manifold. Then, there is no harmonic quasi-isometric map $h:  \Sigma \to Y$.
\el
To prove lemma \ref{lem cusp}, we will use the following result, which will also be crucial for the proof of uniqueness in Theorem \ref{th main}.
\bl\label{lem distance subharmonic}{ \rm \cite[Lemma 5.16]{BH18}}
Let $M$, $Y$ be Riemannian manifolds, and assume that $Y$ has non positive curvature. Let $h_0,h_1:M\to Y$ be harmonic maps. Then, the distance function $m\in M\to d(h_0(m),h_1(m))\in\RR$ is subharmonic.
\el
\noi{\bf Proof of Lemma \ref{lem cusp}} Assuming that there exists a harmonic quasi-isometric map  $h:  \Sigma \to Y$, we will prove that $h$ takes its values in a bounded domain of $Y$. This is a contradiction since $  \Sigma $ is unbounded. 

Let $\hh_t=\{ {\rm Im}\, z=e^t\}/{<\tau>}\subset   \Sigma $  denote the horocyle at distance $t\geq 0$ from $\hh_0$. 
Pick a point $z_0\in\hh_0$ and let $m_0=h(z_0)\in Y$. Since $h$ is harmonic, Lemma \ref{lem distance subharmonic} ensures that the function $\varphi:z\in   \Sigma \to d(h(z),m_0)\in[0,\infty[$ is subharmonic. Since $h$ is quasi-isometric, there exists a constant $k>0$ such that
$\varphi (z)\leq k (t +1)$ for any $z\in\hh_t$, with $t\geq 0$.

Let now $T>0$  and introduce the harmonic function defined on $  \Sigma $ by $$\eta_T(z)= k+k(T+1)\, e^{-T}\, \Im z\, .$$ By construction, $\varphi\leq\eta_T$ on $\hh_0\cup\hh_T$, hence the maximum principle ensures that $\varphi(z)\leq\eta_T(z)$ holds for any $z\in\hh_t$ with $0\leq t\leq T$.
In other words, $d(h(z),m_0)\leq k+k(T+1)\, e^{-T}\, \Im z$ if $T\geq\Im z$. Letting $T\to\infty$ proves that $h$ takes its values in the ball $B(m_0,k)$, a contradiction since $h$ is quasi-isometric. \eproof

\subsection{A second counterexample}\label{subsec:counterex}
We now give another counter-example where the surface  has no cusp, that is when $\Gamma$ has no parabolic element.

Let $\Sigma$ now denote a compact hyperbolic surface whose boundary is the union of two totally geodesic curves of the same length. We consider a non compact hyperbolic surface $S$ obtained by gluing together infinitely many copies of $\Sigma$ along their boundaries, in such a way that $S$ admits a natural action $\tau$ of $\ZZ$ ``by translation''. The quotient $S/\tau$ is then a compact Riemann surface without boundary. Observe that $S$, being quasi-isometric to $\RR$, is Gromov hyperbolic and that  the injectivity radius of $S$ is bounded below but is not a proper function on $S$. 

Our aim is to prove the following.

\bp\label{prop cexZ} There exist quasi-isometric functions $\varphi:S\to\RR$ that are not within bounded distance from any harmonic function.

Thus, there exist quasi-isometric maps $f:S\to\HH^2$ that are not within bounded distance from any harmonic map. 
\ep
\begin{figure}[h]
	\centering
	\includegraphics[width=10.5cm]{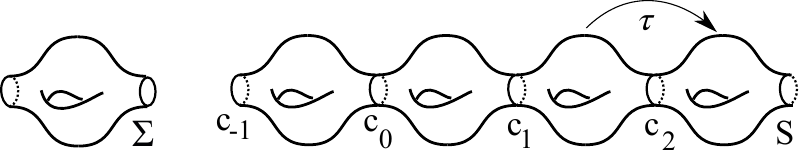}
	\caption*{\small The surface $S$}
\end{figure}
We will first describe all quasi-isometric harmonic functions on $S$.
By applying the afore mentionned theorem by Eeels and Sampson \cite{EellsSampson64} to maps $S/\tau\longrightarrow\RR/\ZZ$ between these compact manifolds, we  construct a harmonic function $\eta_1 : S\longrightarrow\RR$ that satisfies the relation $\eta_1(\tau m) = \eta_1(m)+1$ for any point $m\in S$. We may think of the function $\eta_1$ as a projection from $S$ onto $\RR$. It is a quasi-isometric map. 
\bl Let $\eta :S\longrightarrow\RR$ be a quasi-isometric harmonic function. Then, there exist constants $\a\in\RR^*$ and $\b\in\RR$ such that $\eta =\a \eta_1+\b$.
\el
\proof We denote by $(c_n)_{n\in\ZZ}$ the memories in $S$ of the boundary of $\Sigma$, with $c_n=\tau^n c_0$ ($n\in\ZZ$). All these curves have the same finite length.

By adding a suitable constant to  $\eta _1$, we may assume that $\eta_1$ vanishes at some point  $m_0\in c_0$.  Let $m_n=\tau^n(m_0)\in c_n$ so that $\eta_1(m_n)=n$ for $n\in\ZZ$.

Let $t_n=\eta(m_n)$ ($n\in\ZZ$). Replacing the quasi-isometric function $\eta$ by $-\eta$ if necessary, we may assume that $t_{\pm n}\to \pm \infty$  when $n\to+\infty$. Better, there exist a constant $k>1$ and an integer $N\in\NN$ such that
\begin{equation*}
\begin{gathered}
\hbox{$|\eta_1-n|\leq k$ and $|\eta-t_n|\leq k$} \quad \hbox{on $c_n$ for $n\in\ZZ$,}     \\
\hbox{where $n/k\leq \varepsilon t_{\varepsilon n}\leq kn$ for $\varepsilon =\pm 1$ and $n\geq N$.}
\end{gathered}
\end{equation*}
For any $n\geq 1$, let $(\a_n, \b_n)\in\RR^2$ be the solution of the linear system\\
$$-n\a_n +\b_n=t_{-n} \quad \hbox{and} \quad n\a_n+\b_n=t_n\, .$$
The sequence $(\a_n)$ is bounded since $|\a_n|\leq k$ when $n\geq N$.
The harmonic function $\eta-(\a_n \eta_1+\b_n)$ vanishes at both points $m_n\in c_n$ and $m_{-n}\in c_{-n}$, so that $|\eta-(\a_n \eta_1+\b_n)|\leq k(k+1)$ on $c_n\cup c_{-n}$ if $n\geq N$. By the maximum principle, it follows that $|\eta-(\a_n \eta_1+\b_n)|\leq k(k+1)$ on the compact subset of $S$ cut out by $c_n\cup c_{-n}$.
By applying this estimate at the point $m_0$, we obtain that the sequence $(\beta_n)$ is also bounded. By going to a subsequence, we may thus assume that $\a_n\to\a$ and $\beta_n\to\b$ when $n\to\infty$, so that the limit harmonic function $\eta-(\a \eta_1+\beta):S\to\RR$ is bounded. Since $S$ is a nilpotent cover of a compact Riemannian manifold, a theorem by Lyons-Sullivan \cite[Th.1]{LyonsSull84} ensures that this bounded harmonic function is constant.
\eproof

\vsppp{\bf Proof of Proposition \ref{prop cexZ} \ }
Let $\varphi:S\to\RR$ be a quasi-isometric function such that 

-- $\varphi (m_n)= n$ when $n=\pm 4^{2p}$

-- $\varphi (m_n)= 2n$ when $n=\pm 4^{2p+1}$.

\noi This function $\varphi$ is quasi-isometric, but is not within bounded distance from any function $\a \eta_1+\b$.

We now embed isometrically the real line in the hyperbolic plane as a geodesic $\gamma :\RR\to\HH^2$  and define $f=\gamma\circ \varphi : S\to\HH^2$. Then $f$ is a quasi-isometric map. Assume by contradiction that there exists  a harmonic map $h_1:S\to\HH^2$ within bounded distance from $f$. 
Let $h_2=\sigma\circ h_1:S\to\HH^2$, where $\sigma :\HH^2\to\HH^2$ is the symmetry with respect to geodesic $\gamma$. By Lemma \ref{lem distance subharmonic}, the bounded function $t\in S\to d(h_1(t),h_2(t))\in\plu$ is subharmonic.
Since $S$ is a $\ZZ$-cover of a compact Riemannian manifold, another theorem by Lyons-Sullivan \cite[Th.4]{LyonsSull84} ensures that this bounded subharmonic function is constant.
Therefore, the harmonic map $h_1$ takes its values in a curve which is equidistant to $\gamma$, hence in $\gamma$.
This means that $h_1$ reads as $h_1=\gamma\circ \eta$ with $\eta:S\to\RR$ harmonic within bounded distance from $\varphi$, a contradiction.
\eproof

\section{Proof of Theorem  \ref{th main}}\label{sec:proof}
\begin{quote}
	In this section, we prove  Theorem  \ref{th main}, taking Proposition \ref{prop rho} below for granted.
\end{quote}

Recall that both $X$ and $Y$ are pinched Hadamard manifolds,  that $\Gamma\subset \is (X)$ is a torsion-free convex cocompact discrete subgroup of the group of isometries of $X$ which is not cocompact, and that we let $ M=\Gamma\backslash X$.

Let $f: M\longrightarrow Y$ be a  quasi-isometric map or, more generally, a map that satisfies the hypotheses in Theorem \ref{th main}. 
We want to prove that there exists a unique harmonic map $h: M\longrightarrow Y$  within  bounded distance from $f$.

\subsection{Smoothing the map $f$}\label{subsec:smooth f}
	We first observe that we can assume that the initial map $f: M\to Y$ is  smooth, with bounded covariant derivatives.
\bl\label{lem smooth f} Let $f: M\longrightarrow Y$ be a rough Lipschitz map, namely that satisfies
$$
d(f(m),f(p))\leq c_0\, d(m,p)+c_0
$$ for some constant $c_0$, and any pair of points $m,p\in M$. Then, there exists a $C^\infty$ map $F: M\longrightarrow Y$   with uniformly bounded first and second covariant derivatives, and such that $d(f,F)<\infty$. 
\el
\proof Lift  $f: M\to Y$ to  $\tilde f :X\to Y$. Since $\tilde f$ is still rough-Lipschitz, the construction in \cite[Section 3.2]{BH15} provides a smooth map $\tilde F:X\to Y$ with bounded first and second covariant derivatives, and   within bounded distance from $\tilde f$. This construction being $\Gamma$-invariant, the map $\tilde F$ goes to the quotient and yields the smooth map $F: M\to Y$ we were looking for. \eproof

\subsection{Smoothing the convex core}\label{subsec:smooth core}
	Our goal in this paragraph is to construct the family $(V_R)$  of compact convex neighbourhoods with smooth boundaries of the convex core, on which we solve in  Lemma \ref{lem hR} the bounded Dirichlet problems with boundary value $f$.
\bp \label{prop smooth core} There exists a compact convex set $V\subset  M $ with smooth boundary which is a  neighbourhood of the convex core $K\subset  M $. For any $R>0$, the $R$-neighbourhood $V_R$ of $V$ is also a convex subset of $ M $ with smooth boundary. 
\ep
The proof will rely on Proposition \ref{prop greenwu}, which is due to Greene-Wu.
\bd{\small\bf Strictly convex functions}

\vsp Let $\varphi :M_0\to\RR$ be a continuous function defined on a Riemannian manifold $M_0$.
We say that the function $\varphi$ is strictly convex  if, for every compact subset $L\subset M_0$, there exists a constant $\a>0$   such that, for any unit speed geodesic $t\to c_t\in L$, the function $t\to\varphi (c_t)-\a\,  t^2$ is convex.

\vspace{0.4em}
When $\varphi$ is $C^2$, this definition means  that $D^2\varphi>0$ on $M_0$.
\ed

\bp{\rm\cite[Th.2]{GreeneWu76}}\label{prop greenwu}
Let $M_0$ be a (possibly non complete)  Riemannian manifold and $\varphi :M_0\to\RR$ be a strictly convex function.
Then, there exists a sequence $(\varphi_n)$ of smooth strictly convex functions on $M_0$ that converges uniformly to $\varphi$ on compact subsets of $M_0$.
\ep
The main tool used in \cite{GreeneWu76} to prove Proposition \ref{prop greenwu} is a smoothing procedure called Riemannian convolution.

\bl{\bf\small Strict convexity of $\varphi_{_C}$} \label{lem phiC}

\vsp 
Let $C$ be a non empty closed convex subset of the pinched Hadamard manifold $X$. Then the function $\varphi_{_C} =d^2(.,C)$,  square of the distance function to $C$,  is strictly convex on the complement $M_0=X\smallsetminus C$ of the convex set $C$.
\el
\noi{\bf Proof of Lemma \ref{lem phiC}} We only need to prove that, for every $\varepsilon >0$, there exists an $\a>0$ such that for any unit speed geodesic segment $t\to c_t$ with $d(c_t,C)\geq\varepsilon$, one has
\begin{align*}
2\, \varphi_{_C}(c_{(s+t)/2})- \varphi_{_C}(c_t)- \varphi_{_C}(c_s)\leq -\alpha  (t-s)^2\, .
\end{align*}
Denote by $\pi$ the projection on the closed convex set $C$. Applying the Reshetnyak comparison lemma \ref{lem resh} to the quadrilateral $[c_t,\pi (c_t), \pi (c_s), c_s]$, we are reduced to the well-known case where $C$ is a geodesic segment in the hyperbolic plane $\HH^2$.
\eproof

\bre\label{rem psiC sousharm}   The function $\psi_{_C}=d(.,C)$, distance function to the  convex subset $C$ of $X$, is  convex. Moreover, there exists $\a_1>0$ such that $\Delta \psi_C \geq\a_1$ outside the $1$-neighbourhood of $C$.
This follows from the same arguments as above.\ere

\noi{\bf Proof of Proposition \ref{prop smooth core}}
Applying Lemma \ref{lem phiC} to the convex hull $C:=\conv (\Lambda_\Gamma)$ of the limit set, we obtain that the function $\varphi_{_C}$ is strictly convex on $X\smallsetminus C$. Hence, the function $\varphi_K=d^2(.,K)$ is strictly convex on $M\smallsetminus K$.

We may thus apply Proposition \ref{prop greenwu}  on the manifold $M_0=M\smallsetminus K$ to the function $\varphi_K$, and  obtain a smooth strictly convex function $\varphi_n$ on $M_0$ such that, for every $m$ with $1/2\leq d(m,K)\leq 2$, one has $|\varphi_n(m)-\varphi (m)|\leq 1/2$.

We then  define $V$ as the set $V= K \cup \varphi_n^{-1} ([0,1])$. By construction, $V$ is a convex neighbourhood of the convex core $K$.  Since $\varphi_n$ does not reach a minimum on the boundary $\partial V=\varphi_n^{-1}(1)$, the differential of this convex function  does not vanish on  $\partial V$, so that $V$ has smooth boundary.
\eproof

\subsection{Existence}\label{subsec:existence}

To prove Theorem \ref{th main}, it follows from Paragraph \ref{subsec:smooth f} that we may assume that the  map $f: M\to Y$ we are starting with is not only rough Lipschitz but satisfies,  as well as Condition (b) of Definition \ref{def loc qi},  the stronger condition~: 

\vspp (a$'$) {\sl There exists a constant $c>1$ such that
	\begin{equation}
	\label{eq f smooth}
	\hbox{ $f: M\to Y$ is smooth with $\| Df\|\leq c$ and $\| D^2f\|\leq c$.}
	\end{equation}}
Recall that $V$ is the compact convex neighbourhood, with smooth boundary, of the convex core $K\subset M$ that we constructed in Proposition \ref{prop smooth core}.
\bl\label{lem hR}
For $R>0$, let $V_R\subset M$ denote the $R$-neighbourhood of $V$. Then, there exists a unique harmonic map $h_R:V_R\to Y$ solution of the Dirichlet problem 
\begin{equation*}
\label{eq dirichlet}\hbox{$h_R=f$ on the boundary $\partial V_R$.}
\end{equation*}

\el
\proof This is a consequence of a theorem by R. Schoen \cite[(12.11)]{EeLem78}, since $V_R\subset M$ is a compact manifold with smooth boundary (Lemma \ref{prop smooth core}), $Y$ is a Hadamard manifold and $f: M\to Y$ is a smooth map.
\eproof

\espace
The crucial step in the proof of existence in Theorem \ref{th main} consists in the following uniform estimate.

\bp\label{prop rho}
There exists a constant $\bar\rho>0$ such that $d(f,h_R)\leq\bar\rho$ for any large radius $R> 0$.
\ep
Since the map $f$ is $c$-Lipschitz, we infer that the smooth harmonic maps $h_R$ are locally uniformly bounded. The following statement due to Cheng,  and which is local in nature, provides a uniform bound for their differentials as well.
\bp {\bf\small Cheng Lemma \rm  \cite{Cheng80}}  
Let $Z$ be a $k$-dimensional complete Riemannian manifold with sectional curvature $-b^2\leq K_Z\leq 0$, and $Y$ be a Hadamard manifold. Let 
$z\in Z$,  $r>0$
and let $h:B(z,r)\subset Z\to Y$ be a harmonic $C^\infty$ map 
such that 
the  image 
$ h(B(z,r))$ lies in a ball of radius $r'$. 
Then one has the bound
\begin{equation*}
\|D_zh\|\leq 2^5\,k\,\tfrac{1+br}{r}\, r'\; .
\end{equation*}
\ep
\bc \label{cor cheng} There exists a constant $\kappa \geq 1$ that depends only on $ M$ and $Y$, and with the following property.
Let $R>0$ and assume that $d(f,h_R)\leq\rho$ for some constant $\rho\geq c$. Let $p\in M$ such that $B(p,1)\subset V_R$. Then,
\begin{equation*}
\label{eqncheng}
\|D_ph_R\|\leq \kappa\, \rho\, . 
\end{equation*}
\ec
\proof The map $f$ being $c$-Lipschitz, it follows from $d(f,h_R)\leq\rho$ that 
$$h_R(B(p,1))\subset B(h_R(p),2\rho+c)\subset B(h_R(p),3\rho)\, .$$
Thus, the Cheng Lemma \ref{cor cheng} applies and yields that $\| D_ph_R\|\leq \kappa\rho$, with $\kappa =2^53k(1+b)$. 
\eproof

\vsppp
\noi{\bf Proof of Theorem \ref{th main}}  Applying the Ascoli-Arzela's theorem, it follows from Propositions \ref{prop rho} and \ref{cor cheng} that we can find a increasing sequence of radii $R_n\to\infty$ such that the sequence of harmonic maps $(h_{R_n})$ converges locally uniformly towards a continuous map $h: M\to Y$ which is within bounded distance $\bar\rho$ from $f$.
The Schauder estimates then provide  a uniform bound for the $C^{2,\a}$-norms of the $(h_{R_n})$, hence we may assume that this sequence converges in the $C^2$ norm, so that the limit map $h$ is smooth and harmonic.
We refer to \cite[Section 3.3]{BH18} for more details.
\eproof

\subsection{Boundary estimates}\label{subsec:boundary estimates}
	In this paragraph, we make the first step towards the proof of Proposition \ref{prop rho} by proving  the so-called boundary estimates.

\vspace{0.2em}
The  boundary estimates state that, if the distance $d(h_R,f)$ is very large, this distance is reached at a point which is far away from the boundary $\partial V_R$ of the domain where $h_R$ is defined. More specifically, we have the following.  
\bp \label{prop estimee dist bord}There exists a constant $B$ such that, for any $R\geq 2$ and any point $m\in V_R$, then
$$
d(f(m),h_R(m))\leq B\, d(m,\partial V_R)\, .
$$
\ep
The proof of this proposition is similar to that of \cite[Proposition 3.7]{BH18}, but we must first construct a strictly subharmonic function $\Psi: M\to\plu$ that coincides, outside the $1$-neighbourhood $V_1$ of $V$,  with the distance function $\psi_V=d(.,V)$.  

\bl There exists a constant $\a_1>0$ such that the function $\psi_V=d(.,V)$ satisfies the inequality $\Delta_m\psi_V\geq \a_1$  at any point $m\in M\smallsetminus V_1$ . 
\el
\proof  This follows from Remark \ref{rem psiC sousharm} applied to the convex set $\tilde V\subset X$, which is the lift of $V\subset M$. \eproof

\bc \label{cor Psi}There exists a  continuous function $\Psi: M\to\plu$ which is uniformly strictly subharmonic, namely with $\Delta\Psi\geq\varepsilon$ (weakly) on $ M$ for some $\varepsilon >0$, and such that
$\Psi=\psi_V$  outside $ V_1$.
\ec
\proof If $\nu$ denotes the outgoing unit normal to $V_1$, we have $\psi_V=1$  and $\frac{\partial\psi_V}{\partial \nu}=1$ on $\partial V_1$.
Let $\eta :V_1\to\RR$ be the solution of the Dirichlet problem 
$$
\hbox{$\Delta \eta=1$ on $V_1$, and $\eta=0$ on $\partial V_1$.}
$$
Note that  there exists a constant $c_\eta>0$ such that $-c_\eta\leq \eta\leq 0$ on $V_1$, and $0\leq \frac{\partial \eta}{\partial \nu}\leq c_\eta$ on $\partial V_1$.
Let $\psi_1=1+\frac \eta{c_\eta}:V_1\to\plu$, and define a function $\Psi : M\to\plu$ by letting
$$
\hbox{$\Psi = \psi_1$ on $V_1$ and $\Psi= \psi_V$ outside $V_1$.}
$$
The  function $\Psi$ is positive and continuous on the whole $ M$. Moreover, since $\frac{\partial\psi_V}{\partial \nu}=1\geq \frac{\partial\psi_1}{\partial \nu}$ on $\partial V_1$, it follows that 
$\Delta\Psi\geq \inf (\alpha_1,1/c_\eta)$ weakly on $ M$.
\eproof

\espace\noi{\bf Proof of Proposition \ref{prop estimee dist bord}} 
Let $m\in V_R$. Choose a point $y\in Y$ such that $f(m)$ lies on the geodesic segment $[h_R(m),y]$, and such that $d(y,f(V_R))\geq 1$. Introduce, for some constant $B>0$ to be chosen later on,  the function
$$v:p\in V_R\longrightarrow d(y,h_R(p))-d(y,f(p))-\, (R-\Psi (p))B/2\in \RR\, .
$$
 The choice of $y$ ensures that $d(h_R(m),f(m))=d(y,h_R(m))-d(y,f(m))$.  If the point $m$ lies in $V_R\smallsetminus V_1$, we have $R-\Psi(m)=d(m,\partial V_R)$ while, if $m\in V_1$, the inequality
$R-\Psi (m)\leq R\leq 2(R-1)\leq 2\, d(m,\partial V_R)$ holds. Therefore, we only have to prove that $v(m)\leq 0$.
Since the function $v$ vanishes on the boundary $\partial V_R$, we will be done if we prove that, for a suitable choice of the constant $B$, the function $v$ is subharmonic on $V_R$.

Since $h_R$ is a harmonic map, the function $p\to d(y,h_R(p))$ is subharmonic (Lemma \ref{lem distance subharmonic}). Since $f$ is smooth with uniformly bounded first and second order covariant derivatives, and we chosed $y\in Y$ with $d(y,f(V_R))\geq 1$, it follows that there exists a constant $\beta$ such that the absolute value of the Laplacian of the function $p\in V_R\to d(y,f(p))\in\RR$ is bounded by $\b$ (see \cite[(2.3)]{BH18}).  We infer from Corollary \ref{cor Psi} that
$$
\Delta v\geq 0-\b+\varepsilon B/2 >0
$$ 
if $B$ is large enough, hence the result.\eproof

\subsection {Uniqueness} 
Before going into the main part of  the proof of Proposition \ref{prop rho}, we settle the matter of uniqueness. This is where the non compactness of $M$ is needed. 

\vspace{0.6em}
\noi{\small\bf Proof of uniqueness in Theorem \ref{th main}}

\vsp We rely on arguments in \cite[Section 5]{BH17}. See also \cite[Lemma 2.2]{LiWang98}. Let  $h_0,h_1: M \to Y$ be two  harmonic maps within bounded distance from $f$. We assume by contradiction that  $\delta := d(h_0,h_1)>0$. 

Assume first that we are in the easy case where the subharmonic function $m\in M\to d( h_0(m),h_1(m))\in\plu$ achieves its maximum, hence is constant. As in \cite[Corollary 5.19]{BH17}, it follows  that both maps $ h_0$ and $ h_1$ take their values in the same geodesic of $Y$. Hence each end of $ M$ is quasi-isometric to a geodesic ray. This is a contradiction, since $\Gamma$ being convex cocompact implies that the injectivity radius is a proper function on $ M$.

\vspace{0.4em}
Assume now that there exists a sequence of points $(m_i)_{i\in\NN}$ in  $ M$, that goes to infinity, and such that $d( h_0(m_i),h_1(m_i))\to \delta$. We may also assume that the sequence $(m_i)$ converges to a point $\xi\in\partial_\infty M$.

Since the injectivity radius is a proper function on $M$, it follows from Hypothesis (b) on $f$ that there exist a sequence of radii $r_i\to \infty$ and a constant $c_\xi$   such that $r_i<\inj (m_i)$  and  the restriction of $f$ to each ball $B(m_i,r_i)$ is a quasi-isometric map for some constant $c_\xi$.

Applying \cite[Lemma 5.16]{BH17} ensures that, going if necessary to a subsequence, there also exist two limit $C^{2}$ pointed Hadamard manifolds $(X_\infty,x_\infty)$ and $(Y_\infty, y_\infty)$ with $C^1$ Riemannian metrics  such that the maps
$h_0,h_1 :B(m_i,r_i)\to Y$ respectively converge to quasi-isometric harmonic maps $h_{0,\infty}, h_{1,\infty}:X_\infty\to Y_\infty$ with $d(h_{0,\infty}(x),h_{1,\infty}(x))=\delta$ for every $x\in X_\infty$. 

Applying again \cite[Corollary 5.19]{BH17}, we infer that both maps $h_{0,\infty}, h_{1,\infty}$ take their values in the same geodesic of $Y$.  This is a contradiction, since both maps $h_{i,\infty}$ are quasi-isometric.
\eproof

\section{Lower bound for harmonic measures}\label{sec:harmonic measure, lower}
\begin{quote}
	To prove Proposition \ref{prop rho} in Section \ref{sec:interior estimates}, we will need uniform bounds for the harmonic measures on specific domains of $ M$.  In this section, we deal with the lower bounds.
\end{quote}
\subsection{Harmonic measures}\label{subsec: harmonic measure}
Assume that $M$ is a Riemannian manifold, and let $W \subset M$ be a relatively compact domain with smooth boundary. For any continuous function $u: \partial W \to\RR$, there exists a unique continuous function $\eta_u:\overline W \to\RR$ which is smooth on $W $, and is solution to the Dirichlet problem
$$
\hbox{$\Delta \eta_u =0$ on $W $ \quad and \quad $\eta_u=u$ on $\partial W $.}
$$
This gives rise to a family of Borel probability measures $\sigma_{m,W } $ supported on $\partial W $, indexed by the points $m\in W $ and such that, for any continuous function $u\in C^0(\partial W )$~:
\begin{equation}\label{eq mesure harmonique}
\eta_u(m)=\int_{\partial  W } u(z)\, d\sigma_{m, W } (z)\,  .
\end{equation}
The measure $\sigma_{m, W }$ is the harmonic measure of $ W $ relative to the point $m$.

\espace
In our previous paper \cite{BH17}, we worked in pinched Hadamard manifolds, and obtained the following uniform upper and lower bounds for the harmonic measures on balls relative to their center.
\bt{\rm\cite{BH17}} \label{th harmonic measure on X} Let $0<a\leq b$ and $k\geq 2$. There exist positive constants $C,s$ depending only on $a,b$ and $k$, and with the following property.

Let $X$ be a $k$-dimensional pinched Hadamard manifold, whose sectional curvature satisfies $-b^2\leq K_X\leq -a^2$.

Then for any point $x\in X$, any radius $R>0$ and  angle $\theta\in [0,\pi/2]$, the harmonic measure $\sigma_{x,R}$ of the ball $B(x,R)$ relative to the center $x$ satisfies
\begin{equation}\label{eq estimees uniformes mesure X}
\frac 1C\,\theta^s\leq \sigma_{x,R} (\C_x^\theta)\leq C\,\theta^{1/s}
\end{equation}
where $\C_x^\theta\subset X$ denotes any cone with vertex $x$ and angle $\theta$.
\et
Since the measure $\sigma_{x,R}$ is supported on the sphere $S(x,R)$, the expression $\sigma_{x,R} (\C_x^\theta)$ means
$\sigma_{x,R} (\C_x^\theta\cap S(x,R))$.

\espace
To prove Theorem \ref{th main}, we will need similar estimates in the quotient manifold $ M$. The lower bound is provided in the next paragraph. The upper bound will require more work and will be carried out in Section \ref{sec:upper bound}.

\subsection{Harmonic measures on $ M$}\label{subsec: harmonic measure quotient}
The description of the positive harmonic functions on the quotient $M=\Gamma\backslash X$ of a pinched Hadamard manifold by a convex cocompact group is due to Anderson-Schoen in \cite[Corollary 8.2]{AndersonSchoen85}. 

Our goal in this paragraph is to obtain the following lower bound for the mass of a ball of fixed radius $\lambda >0$,  with respect to the harmonic measures  on  suitable domains of $ M$. 
\bp\label{prop lower bound harmonic} {\bf\small Lower bound for harmonic measures on $ M=\Gamma\backslash X$}

\vsp Let  $\lambda >0$ and $L\geq 1$. Then, there exists a constant $s(\lambda , L)>0$ with the following property.
For any pair of points $m,q$ in $ M$ with $ 0<d(m,q)\leq L $, there exists a
bounded domain with smooth boundary $ W  _{m,q}$ with $m\in  W  _{m,q}$ and $q\in \partial  W  _{m,q}$, and whose harmonic measure relative to the point $m$ satisfies the inequality
\begin{equation}\label{eq minoration mesure}
\sigma_{m, W  _{m,q}}\bigl(  B(q,\lambda)\bigr) \geq s(\lambda, L)\, .
\end{equation}
\ep
This statement will derive from a compactness argument, together with the following continuity lemma. 

\bl \label{lem dependance mesure harmonique}
Let $  W \subset\RR^k$ be a bounded domain with $C^{2}$ boundary, and $(g_n)_{n\in\NN}$ be a sequence of Riemannian metrics on $\overline  W $ that  converges, in the $ C^2(\overline  W )$ sense, to a Riemannian metric $g_\infty$ on $\overline  W $.

Let $u\in C^0(\partial  W )$ be a continuous function on the boundary. For each $n\in\NN\cup\{\infty\}$, denote by $\eta_n\in C^2(  W )\cap C^0(\overline  W )$
the solution of the Dirichlet problem with fixed boundary value $u$ for the metric $g_n$, namely such that
$$
\Delta_n \eta_n=0 \qquad (\eta_n)_{|\partial  W }=u\, .
$$
Here $\Delta_n$ denotes the Laplace operator corresponding to the metric $g_n$.

Then, the sequence $(\eta_n)_{n\in\NN}$ converges uniformly on $\overline  W $ to $\eta_\infty$.
\el
\proof Introduce the solution $\greta :  W \to\RR$ of the boundary value problem
$$
\Delta_\infty \greta=1 \qquad \greta_{|\partial  W }=0\, .
$$
 Let $\varepsilon >0$. If $n$ is large enough, we have
$$
\hbox{$\Delta_n (   \eta_\infty-\eta_n +\varepsilon\greta)\geq 0$ and $ \Delta_n (   \eta_\infty-\eta_n -\varepsilon\greta)\leq 0$.}
$$
The claim follows. Indeed, both functions $\eta_\infty-\eta_n \pm\varepsilon\greta$ vahishing on the boundary $\partial  W $,
 the maximum principle ensures that
\begin{equation}
|\eta_\infty-\eta_n|\leq \varepsilon\, \sup |\greta|\, .  \plouf
\end{equation}

\espace
Each domain $W_{m,q}\subset M$ will either be a ball, or  the image under a suitable diffeomorphism of a fixed domain  $  W \subset\RR^k$. This diffeomorphism will be defined using  the normal exponential map along a geodesic segment containing $[m,q]$. We first observe the following, where   $\inj (M) >0$ denotes the injectivity radius of $M$.
\bl\label{lem geo injective}
(1) Let $[m,q]\subset M$ be a minimizing geodesic segment. 
Then, the  extended geodesic segment $I_{mq}=[m_q, q]$ defined by the conditions $[m,q]\subset [m_q, q]$ and  $d(m,m_q)=\inj (M)/2$ is still injective.

\vspace{0.3em}
(2) For any compact subset $Z\subset M$, there exists $0<r\leq 1$ such that, when $[m,q]$ is a minimizing geodesic segment with $m\in Z$ and  $d(m,q)\leq L$, the normal exponential map $\nu_{mq}$ along $I_{mq}$ is a diffeomorphism from the bundle of normal vectors to $I_{mq}$ with norm at most $r$  onto its image.
\el
\proof (1) derives easily from the definition of the injectivity radius.

(2)  follows since the $L$-neighbourhood of $Z$ is also compact.
\eproof

\espace
Let $\a =\inj (M)/2L$ and introduce the segment $J=[-\a,1]\subset\RR$. We choose $W$ to be  a convex domain of revolution $  W \subset J\times B(0,r)\subset\RR\times\RR^{k-1}$ whose boundary $\partial   W $ is smooth and contains both points $(-\a,0)$ and $(1,0)$.

\espace
\noi{\bf Proof of Proposition \ref{prop lower bound harmonic} }
 We may assume that $\lambda\leq\inj (M)$. 
 
  If $d(m,q)\leq \lambda /2$, then choose $W_{m,q}$ to be the ball with center $m$ and radius $d(m,q)$, so that $W_{m,q}\subset B(q,\lambda)$.

 We now assume that  $d(m,q)>\lambda /2$. Since the injectivity radius is a proper function on $ M$, the set
$$
Z=\{ m\in  M \, |\; \inj (m)\leq L +\lambda\} 
$$
is a compact subset of $ M$. 

\vspace{0.5em} 
Assume first that the point $m$ does not belong to the compact set $Z$. Then, the  ball $B(m,L+\lambda)$ is isometric to a ball with radius $L+\lambda$ in the Hadamard manifold $X$ while $B(q,\lambda)\subset B(m,L+\lambda)$. Choosing $ W  _{m,q}=B(m,d(m,q))$, the required estimate follows easily from the lower bound in Theorem \ref{th harmonic measure on X}.

Assume now that $m\in Z$, and that $\lambda/2\leq d(m,q) \leq L$. Pick a minimizing geodesic segment $[m,q]$. Identify $\RR$ with the geodesic line containing $[m,q]$ through the constant speed parameterization $c:\RR\to M$ defined by $c(0)=m$ and $c(1)=q$, so that $c(J)\subset I_{mq}$. Introduce
$W_{m,q}=\nu_{m,q }(W)$. By construction, $W_{m,q}$  is a bounded domain of $M$ with smooth boundary, such that $m\in W_{m,q}$ and $q\in \partial W_{m,q}$.

 Let us prove that the harmonic measures $\sigma_{m, W  _{m,q}}\bigl(  B(q,\lambda)\bigr) $ are uniformly bounded below.
We proceed by contradiction and assume that there exist two sequences of points $m_n\in Z$, and $q_n\in  M$ with $\lambda /2 \leq d(m_n,q_n)\leq L$, and such that $\sigma_{m_n, W  _{m_n ,q_n}}\bigl(  B(q_n,\lambda)\bigr) \longrightarrow 0$ when $n\to\infty$.

Since  $Z$ is compact, we may assume that $m_n\to m_\infty\in Z$, that $q_n\to q_\infty$ with $m_\infty\neq q_\infty$, and that the sequence of minimizing geodesic segments $([m_n,q_n])$  converges to a minimizing geodesic segment $[m_\infty, q_\infty]$. Denoting by $g_n$ ($n\in\NN\cup\{\infty \}$) the Riemannian metrics on $W$ obtained by pull-back of the Riemannian metric of $ M$ under the map $\nu_{m_n,q_n }$, we may even assume that $g_n\to g_\infty$ in the $C^2$ sense on $\overline W$. Hence Lemma \ref{lem dependance mesure harmonique} yields $\sigma_{m_\infty, W  _{m_\infty,q_\infty}}\bigl(  B(q_\infty,\lambda)\bigr)=0$,  a contradiction to the maximum principle.
\eproof
\begin{figure}[h]
	\hspace{1.8em}\includegraphics[width=11cm]{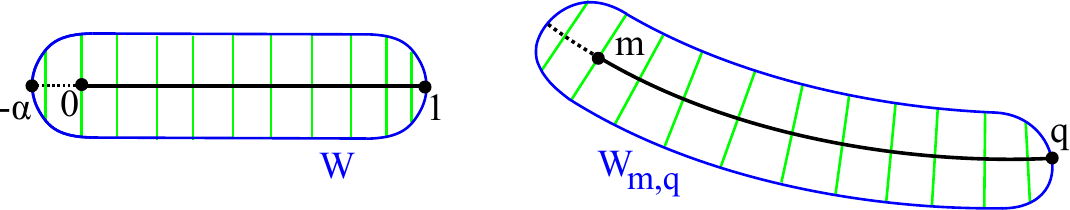}
	\caption*{\small Construction of the domain $W_{m,q}$}
\end{figure}

\section{Upper bound for harmonic measures}\label{sec:upper bound}
\begin{quotation}
	The main goal of this section is to obtain, in Corollary \ref{cor upper bound harmonic},  the upper bound for the harmonic measures on $ M$ needed for the proof of Proposition \ref{prop rho}.
\end{quotation}
One of the major technical tools for estimating or constructing harmonic functions on Hadamard manifolds is the so-called Anderson-Schoen barriers. Given two opposite geodesic rays in a pinched Hadamard manifold, the corresponding Anderson-Schoen barrier is a positive superharmonic function  that decreases exponentially along one of the geodesics rays, and is greater than 1 on a cone centered around the other geodesic ray  \cite{AndersonSchoen85}.

Our first step is to obtain an analogous to the Anderson-Schoen barrier functions for the quotient manifold $ M=\Gamma\backslash X$ in Proposition \ref{prop mesure anderson schoen}. 
This will rely on the work by Ancona \cite{Ancona87}, Anderson \cite{Anderson83}  and Anderson-Schoen \cite{AndersonSchoen85}.

 In the whole section, $X$ will denote our pinched Hadamard manifold satisfying \eqref{eq pinching}. The torsion-free convex cocompact subgroup $\Gamma$ of $\is (X)$ will only come into the picture 
starting from Paragraph \ref{subsec:embedded half}.

\subsection{Harmonic measures at infinity on $X$}\label{subsec: omega sur X}
	We first recall some fundamental, and by now classical, results concerning harmonic measures on Hadamard manifolds. 

\vspace{0.3em}
In the sequel $T,\a,\ci,c_r\geq 1$ will be various constants such that $T,\a,\ci$ depend only on the pinched Hadamard manifold $X$, and $c_r$ also depends on the distance $r>0$.

\espace
Let us first recall the following Harnack-type  inequality, due to Yau.
\bl \label{lem harnack}{\rm\cite{Yau75}}  For a positive harmonic function 
$\eta :B(x,1)\subset X\to \mathopen] 0,\infty\mathclose[$ defined on a ball with radius $1$, one has
$|D_x (\log\eta)|\leq c_1$.
\el
Another fundamental tool is the Green function $G:X\times X\to\mathopen]0,\infty]$. This Green function is continuous on $X\times X$, and is uniquely defined by the conditions
\begin{align*}
\Delta_y G(x,y)&=-\delta_x\cr
\lim_{y\to \infty}G(x,y)&=0\, 
\end{align*}
for every $x\in X$.  One can prove that $G$ is symmetric i.e. $G(x,y)=G(y,x)$ for  $x,y$ in $X$. 
Moreover, the Green function satisfies the following estimates.
\bp\label{prop ppte de G} 
1.   If $d(x,y)\geq 1$, one has
\begin{align}\label{eq Log G qi}
\ci^{-1}\,d(x,y)-\ci\leq \log( 1/G(x,y))\leq \ci\, d(x,y)+\ci\, .
\end{align}
2.  Let  $x,y,z\in X$ such that $d(x,y)\geq 1$ and $d(y,z)\geq 1$. One has
\begin{align}\label{eq log G distance mino}
\ci^{-1}\, G(x,y)\, G(y,z)\leq G(x,z) \, .
\end{align}
If  $d(x,z)\geq 1$ and $d(y,[x,z])\leq r$, one has
\begin{align}\label{eq log G distance majo}
G(x,z)\leq c_r\, G(x,y)\, G(y,z)\, .
\end{align}
\ep
The first assertion is (2.4) in Anderson-Schoen  \cite{AndersonSchoen85}. 
The  second assertion follows from \eqref{eq Log G qi}, using Harnack inequality and the maximum principle, while the third one is  Ancona's inequality  \cite[Theorem 5]{Ancona87}.

\espace Let us now turn to the bounded harmonic functions on $X$. Anderson proved in \cite{Anderson83} that the Dirichlet problem at infinity on $X$ has a unique solution for any continuous boundary value. Hence there exists, for every point $x\in X$, a unique Borel measure $\omega_x$ on $\partial_\infty X$ such that, for any continuous function $u\in C^0(\partial_\infty X)$, the function
$$
\eta_u :x\in X\to \int_{\partial_\infty X} u(\xi)\, d\omega_x(\xi)\in\RR
$$
is harmonic on $X$ and extends continuously to $\overline X$ with boundary value at infinity equal to $u$. The measure $\omega_x$ is the harmonic measure on $\partial_\infty X$ at the point $x$.
The Harnack inequality of Lemma \ref{lem harnack} ensures that two such harmonic measures $\omega_x$ and $\omega_y$ are absolutely continuous with respect to each other, and that their Radon-Nikodym derivatives $\frac{d\omega_y}{d\omega_x}$ are uniformly bounded when $d(x,y)\leq 1$. We will also need a control on these Radon-Nikodym derivatives for $d(x,y)\geq 1$, that will be given in Lemma \ref{lem dmgo}.

\espace
In \cite{AndersonSchoen85}, Anderson-Schoen study the positive harmonic functions on $X$, and provide an identification of the Martin boundary of $X$ with the boundary at infinity $\partial_\infty X$. More precisely, they obtain the following results. 

Fix a base point $o\in X$ and 
introduce the normalized Green function at the point $o\in X$ with pole at $z\in X$, which is defined by
\begin{equation*}
k(o,x,z)=\frac{G(z,x)}{G(z,o)}\, .
\end{equation*}
Letting the point $z\in X$ converge to $\xi\in\partial_\infty X$, the limit
\begin{equation*}
k(o,x,\xi)=\lim_{z\to\xi}k(o,x,z)\, 
\end{equation*}
exists and $x\to k(o,x,\xi)$ is now a positive harmonic function on the whole $X$ that extends continuously to the zero function on $\partial_\infty X\smallsetminus\{\xi\}$, and such that $k(o,o,\xi)=1$.  

In \cite[Theorem 6.5]{AndersonSchoen85} Anderson-Schoen prove that these positive harmonic functions are the minimal ones. Using the Choquet representation theorem, they  provide the following Martin representation formula. For any positive harmonic function $\eta$ on $X$, there exists a unique finite positive Borel measure $\mu_\eta$ on $\partial_\infty X$ such that, for every $x\in X$,
$$
\eta (x)=\int_{\partial_\infty X} k(o,x,\xi)\, d\mu_\eta (\xi)\, .
$$
The minimal harmonic functions relate to the harmonic measures at infinity.
\bl\label{lem dmgo}
1. Let $o,x\in X$. The following holds for $\omega_o$-a.e. $\xi\in \pX$~:
\begin{equation}\label{eq rad der k}
k(o,x,\xi)=\frac{d\omega_x}{d\omega_o}(\xi)\, .
\end{equation}
2. For  $o,x\in X$ and $\xi\in\pX$ with $d(o,x)\geq 1$, we have
$$
k(o,x,\xi)\, G(o,x)\leq \ci\, .
$$
If moreover  $d(x,[o,\xi[)\leq r$, then we also have
$$
c_r^{-1}\leq k(o,x,\xi)\, G(o,x)\, .
$$
\el
\proof 1. is proved in \cite[\S 6]{AndersonSchoen85}.

\noi 2. follows readily from  \eqref{eq log G distance mino} and \eqref{eq log G distance majo}, with $z\in [o,\xi[$ and letting $z\to\xi$.\eproof

\espace\noindent
Lemma \ref{lem dmgo} asserts that, if $d(o,x)\geq 1$ then, for every  $\xi\in\pX$ which is in the shadow of the ball $B(x,r)$ seen from $o$, all the densities $\frac{d\omega_x}{d\omega_o}(\xi)$ are close to $(G(o,x))^{-1}$ hence do not depend too much on $\xi$.

\subsection{The action of $\is (X)$ on $\pX$}
	We now investigate, using Lemma \ref{lem dmgo}, the action of $\is (X)$ on the harmonic measures at inifinity. 

\vspace{0.3em}
Introduce the function defined, for any pair of points $x,y\in X$, by $$d_G(x,y)=\log^+(1/G(x,y))$$ where $\log^+(t)=\sup (\log t,0)$ denotes the positive part of the logarithm.  Proposition \ref{prop ppte de G} tells us that, at large scale, the function $d_G$ behaves roughly like a distance that would be quasi-isometric to the Riemannian distance $d$.  In particular,  \eqref{eq log G distance mino} ensures that there exists a constant $\ci$ such that the following weak triangle inequality holds for every $x,y,z\in X$:
\begin{equation}\label{eq inegalite triang dG}
d_G(x,z)\leq d_G(x,y)+d_G(y,z)+\ci\, .
\end{equation}
Although $d_G$ is not exactly a distance on $X$, we will thus nevertheless agree to think of $d_G$ as of the Green distance.
We would like to mention that Blachère-Haïssinski-Mathieu  \cite{BHaissM11} already used such a Green distance  in the similar context of random walks on hyperbolic groups.

\espace
From now on, we choose a base point $o\in X$. We associate to $d_G$ an analog to the Busemann functions by letting, for $\xi\in\pX$ and $x\in X$~:
\begin{equation*}
\beta_\xi (o,x)=-\log k(o,x,\xi)=\lim_{z\to\xi} d_G(x,z)-d_G(o,z)\, .
\end{equation*}
These Busemann functions relate to the harmonic measures at infinity, since \eqref{eq rad der k} reads as
\begin{equation}\label{eq rad der beta}
\frac{d\omega_x}{d\omega_o}(\xi)=e^{-\beta_\xi(o,x)}\, . 
\end{equation}
Define the length of $g\in\is (X)$ as $|g|_o= d_G(o,go)$. We want to compare $\beta_\xi(o,g^{-1}o)$ with $|g|_o$.
\bno When $g\in\is (X)$ does not fix the point $ o$, we introduce the endpoints $\xi_{_+}\!(g),\xi_{_-}\!(g)\in\pX$ of the geodesic rays with origin 
$o$ that contain respectively the points $go$ and $g^{-1}o$.
\eno

\bl\label{lem g et beta} (1)  For $g\in\is (X)$ and $\xi\in\pX$, one has $|g|_o=|g^{-1}|_o$ and
\begin{equation*}
|\beta_\xi(o,g^{-1}o))|\leq |g|_o +\ci\, .
\end{equation*}
(2) For every $\varepsilon >0$ there exists a constant $a_\varepsilon \geq 1$ such that, when  $g\in\is (X)$ 
and  $\xi\in\pX$  satisfy $\angle_o(\xi, \xi_{_-}\!(g))\geq\varepsilon$, then
$$
|g|_o\leq \beta_{\xi}(o,g^{-1}o)+a_\varepsilon\, .
$$
\el
\proof (1) We observe that, since  the Green function $G$  is symmetric and invariant under  isometries, we have $d_G(y,x)=d_G(x,y)=d_G(gx,gy)$ for every $g\in\is (X)$ and $x,y\in X$. 
The equality $|g|_o=|g^{-1}|_o$ follows.

Using the weak triangle inequality \eqref{eq inegalite triang dG} for $d_G$ yields
$$
\beta_\xi (o,g^{-1}o) =\lim_{z\to\xi} d_G(g^{-1}o,z)-d_G(o,z)\leq d_G(g^{-1}o,o) +\ci =|g|_o+\ci\, .
$$
The lower bound  follows  by observing that the invariance of $d_G$ under isometries ensures that
$$\beta_{\xi}(o,g^{-1}o)=\beta_{g\xi }(go,o)=-\beta_{g\xi }(o,go)\geq -|g|_o-c_\circ\, .$$ 

\noi\begin{minipage}{7.2cm} 
	(2) We may suppose that 
	$|g|_o>0$ so that $|g|_o=-\log G(o,g^{-1}o)$.
	Since $K_X\leq -a^2$,   the condition  $\angle_o(\xi, \xi_{_-}\!(g))\geq\varepsilon$ ensures that the distance of the point $o$ to the geodesic ray $\mathopen[ g^{-1}(o),\xi\mathclose[$ is bounded above by a constant $r_\varepsilon$ that depends only on $\varepsilon$.  
	 Lemma \ref{lem dmgo} (2) yields 
	$|g|_o\leq \beta_{\xi}(o,g^{-1}o) +\log c_{r_\varepsilon}$.
	\eproof	
	
	\vspace{0.7em}
	The following corollary provides useful estimates for action of isometries on the harmonic measures on $\partial_\infty X$. 
	\end{minipage}
\begin{minipage}{7.5cm}
	\hspace{0.4em}		\includegraphics[width=5cm]{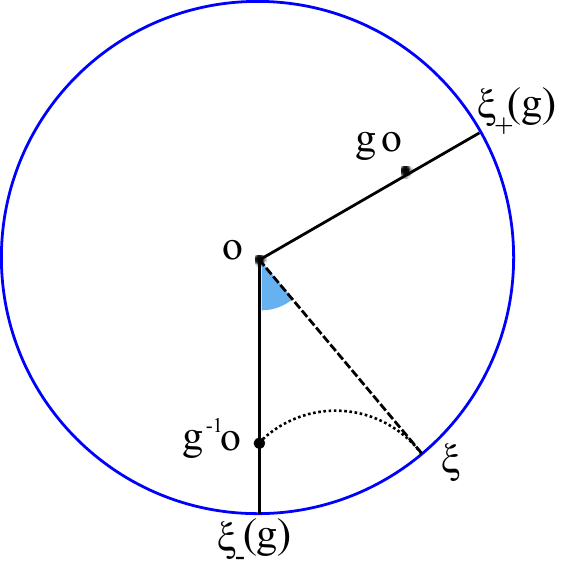}
\end{minipage}

\bc\label{cor omega gammaA}
For $A\subset \pX$ a measurable set and $g\in \is (X)$, one has
$$e^{-c_\circ} e^{-|g|_o}\omega_o(A)\leq \omega_o(g A)\leq e^{c_\circ}  e^{|g|_o}\omega_o(A)\, .$$
If we assume that $\angle_o(\xi, \xi_{_-}\!(g))\geq\varepsilon$ for every $\xi\in A$, one has
$$
\omega_o(g A)\leq e^{a_\varepsilon}  e^{-|g|_o}\omega_o(A)\, 
$$
where $a_\varepsilon$ is the constant in Lemma \ref{lem g et beta}.
\ec
\proof Immediate consequence of Lemma \ref{lem g et beta} and Equality \eqref{eq rad der beta}. \eproof

\subsection{Harmonic measures of cones in $X$}\label{subsec:cones X}
	We now recall estimates for the harmonic measures at infinity, that are due to Anderson-Schoen.

\bd\label{def cone} Let $o\in X$, $\xi\in \partial_\infty X$ and $\theta\in [0,\pi]$. 
The closed cone ${\cal C}_{o\xi}^\theta\subset X$ with vertex $o$, axis  $[o,\xi\mathopen[$ and angle $\theta$ is the union of all the geodesic rays $[o,\zeta\mathopen[$  whose angle with $[o,\xi\mathopen[$ is at most $\theta$. The trace of the cone ${\cal C}_{o\xi}^\theta$ on the sphere at infinity $\pX$ will be denoted by ${\cal S}_{o\xi}^\theta$.
\ed
By analogy with the case where $X$ has constant curvature, we will take the liberty of calling  a cone ${\cal D}_{o\xi}:={\cal C}_{o\xi}^{\pi /2}\subset X$  with angle $\theta =\pi /2$ a closed half-space, with vertex $o$, and its boundary ${\cal H}_{o\xi}$ a hyperplane, also with vertex $o$. We will denote by ${\cal S}_{o\xi}$  the trace of the half-space ${\cal D}_{o\xi}$ on the boundary at infinity, and we will call it a half-sphere at infinity seen from the point $o$.

\vspace{0.3em}
Although half-spaces  in the pinched Hadamard manifold $X$ may not be convex, the following lemma tells us that they are not far from being so.
\bl\label{lem bowditch quasi convex}  {\rm\cite[Prop.2.5.4]{Bowditch95}}
There exists a constant $\lambda$, that depends only on the pinching constants of $X$, such that the convex hull of a half-space  $\D\subset X$  lies within its $\lambda$-neighbourhood~: ${\rm Hull}\, (\D)\subset {\cal V}_\lambda (\D)$.
\el
\proof This statement follows from Proposition \ref{prop bowditch}, due to Bowditch.
Observe indeed that $\D$ is included in the join of its half-sphere ${\cal S}$ at infinity, and that this join lies in the $\delta_X$-neighbourhood of $\D$.
\eproof

\espace
The following uniform  bounds for harmonic measures of cones in pinched Hadamard manifolds are due to Anderson-Schoen.
\bl\label{lem kifer-ledrappier}There exists a constant $c_\circ \geq 1$ such that one has, for any point $o\in X$ and any $\xi\in\pX$,
\begin{equation*}
\omega_o({\cal S}_{o\xi}^{\pi /4})\geq 1/c_\circ\, .
\end{equation*}
\el
A more precise statement is given by Kifer-Ledrappier in \cite[Theorem 4.1]{KiferLedrappier}.
\bno \label{not demi X}We now fix a base point $o\in X$. When $\xi\in\pX$, we denote by $t\in \RR\to x_\xi ^t\in X$ the unit speed geodesic with origin $o$ that converges to $\xi\in\pX$ in the future. We will let $\D_\xi^t$ stand for the half-space ${\cal D}_{x^t_\xi\xi}$ with vertex $x^t_\xi$ and axis $[x^t_\xi,\xi\mathopen[$, and  will denote accordingly by $\H_\xi^t$ and ${\cal S}_\xi^t$ the corresponding hyperplane and half-sphere at infinity.  
\eno
\bp \label{prop a-s k-l} There exist a distance $T>0$ and two constants  $\a >0$  and $c_\circ\geq 1$   such that the following holds for every $\xi\in\partial_\infty X$~:
\begin{alignat*}{2}
\omega_{x^{-t}_\xi}({\cal S}_{o\xi})&\leq c_\circ\, e^{-\a t} & \quad & \text{for every $t\geq 0$}\cr
\omega_y({\cal S}_{o\xi})&\geq 1/c_\circ &       & \text{for every $y\in {\cal D}_\xi^T$.}
\end{alignat*}
\ep 
\proof The first assertion is   \cite[Corollary 4.2]{AndersonSchoen85}.

\vspace{0.2em}
\noi\begin{minipage}{7.cm}
	Thanks to the upper bound on the curvature $K_X$ of $X$, there exists a distance $T>0$ that depends only on $X$ such that the half-space ${\cal D}_\xi^{T}$ is seen from the point $o$ under an angle at most $\pi/4$, namely such that
	${\cal D}_\xi^{T}\subset {\cal C}_{o\xi}^{\pi /4}$. 
	
	If now $y\in {\cal D}_\xi^{T}$ and $\zeta_y\in\pX$ denotes the endpoint of the geodesic ray such that $y\in [o,\zeta_y\mathclose[$, it follows that $\C_{y\zeta_y}^{\pi /4}\subset \D_{o\xi}$. 
	Lemma \ref{lem kifer-ledrappier} yields  $\omega_y({\cal S}_{o\xi})\geq \omega_y({\cal S}_{y\zeta_y}^{\pi /4}) \geq 1/c_\circ $.  \eproof
\end{minipage}
\begin{minipage}{6.8cm}
	\hspace{0.5em}	\includegraphics[width=5.5cm]{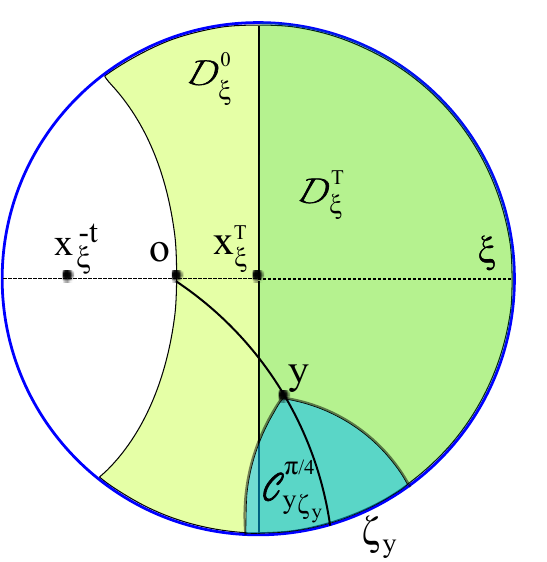}
\end{minipage}

\subsection{Geometry in embedded half-spaces in $ M$}\label{subsec:embedded half}
	In this paragraph, we  introduce embedded half-spaces in the quotient $ M=\Gamma\backslash X$ (Proposition \ref{prop tau0}), that  will be needed in the sequel of this paper.

\vspace{0.5em}
Recall that we fixed a base point $o\in X$. We keep Notation \ref{not demi X}.

\bno
We now introduce  the projection $m_0\in M=\Gamma\backslash X$ of our base point $o\in X$. When $\xi\in \partial_\infty X$, we will denote by $t\in\RR\to m^t_\xi\in M$ the (perhaps non minimizing) geodesic obtained by projection of the geodesic $t\in\RR\to x_\xi^t\in X$.
\eno
\bd\label{def closed half}  
A closed embedded half-space 
with vertex $m_\xi^t$ is the projection in $ M$ of a closed half-space $\D_\xi^t\subset X$ that embeds in $ M$, namely that satisfies  $ \overline\D_\xi^t\cap \gamma  \overline\D_\xi^t=\emptyset$ 
for every non trivial element $\gamma\in \Gamma$, where $\overline\D_\xi^t\subset\overline X$ denotes the closure of $\D_\xi^t$ in the compactification of $X$. 
\ed
We first remark that there exist many embedded half-spaces in $ M$. 
Choose a relatively compact subset $\Omega_0\subset \Omega_\Gamma$ of the domain of discontinuity. 
\bl\label{lem t0 pour embedded}
There exists $t_0>0$   such that, for every $\xi\in \Omega_0$ and every $t\geq t_0$, the half-space $\D_\xi ^t\subset X$ embeds in $ M$.
\el
\proof We proceed by contradiction, and assume that there exist sequences $t_n\to +\infty$, $\xi_n\in\Omega_0$, $z_n\in \overline\D_{\xi_n}^{t_n}$ and $\gamma_n\in\Gamma^*$ such that $\gamma_nz_n\in \overline\D_{\xi_n}^{t_n}$ for every $n\in\NN$. By compactness of $\Omega_0$, we may assume that the sequence $(\xi_n)$ converges to some $\xi_\infty\in\Omega_0$. Since $t_n\to\infty$ and $\xi_n\to\xi_\infty$,  both sequences $(z_n)$ and $(\gamma_nz_n)$ converge to $\xi_\infty$ in $\overline X$. Since the action of $\Gamma$ on $X\cup\Omega_\Gamma$ is properly discontinuous, and $\Gamma$ is torsion-free, it follows that $\gamma_n$ is trivial for $n$ large, a contradiction.
\eproof

\bno\label{not demi} For $\xi\in \Omega_0$ and $t\geq t_0$, we will denote  by the Roman letters $D_\xi^t\subset M$ the closed embedded half-space with vertex $m_\xi^t$ and by $H_\xi^t$ its boundary in $ M$, respectively obtained as the projections of 
$\D_\xi^t $ and  $\H_\xi^t$.
\eno
Our goal in the remaining of this section is to prove   that an embedded half-space in  the quotient $ M=\Gamma\backslash X$ is not far from being geodesically convex~:
\bp \label{prop tau0} There exists $\tau_0\geq 1$ such  that, for every $\xi\in \Omega_0$ and $t\geq t_0$, and 
for every pair of points $p_1,p_2\in D_\xi ^{t+\tau_0}$, there exists only one minimizing geodesic segment $[p_1,p_2]\subset M$, and it  lies in $D_\xi ^t$.
\ep
Proposition \ref{prop tau0} relies on an analogous property for the half-space $\D_\xi^{t+t_0}$ in the pinched Hadamard manifold $X$ that we proved in Lemma \ref{lem bowditch quasi convex}.
We first state  an elementary property of obtuses triangles in $X$. 
\bl\label{lem right angle triangle in X} There exists a constant $\hat\delta_0$ that depends only on $X$ such that if $[y,x,z]\subset X$ is a triangle  with an angle at least $\pi/2$ at the vertex $x$, then  $d(y,z)\geq d(x,y)+d(y,z)-\hat\delta_0$.
\el
\proof Same as for Lemma \ref{lem quad X}.\eproof

\espace
\noi{\bf Proof of Proposition \ref{prop tau0}} Let $\tau_0>0$ and $p_1,p_2\in D_\xi ^{t+\tau_0}$. Consider a minimizing geodesic segment $[p_1,p_2]\subset  M$. Lift the points $p_1,p_2\in D_\xi ^{t+\tau_0}$ as $y_1,y_2\in \D_\xi^{t+\tau_0}$.
Then, there exists $\gamma\in \Gamma$ such that $[p_1,p_2]$ lifts as a geodesic segment $[y_1,\gamma y_2]\subset X$.

\noi
\begin{minipage}{7cm}
	\vspace{0.8em}\noi	Assume first that $\gamma =e$. If we choose $\tau_0\geq\lambda_X$, the $\lambda_X$-neighbourhood of the half-space $\D_\xi^{t+\tau_0}$ lies in $\D_\xi^{t}$. Hence, it follows from Lemma \ref{lem bowditch quasi convex} that $[y_1,y_2]\subset \D_\xi^{t}$ so that $[p_1,p_2]\subset D_\xi ^t$.
	
	\vspace{0.6em}\noi
	Proceed now by contradiction and assume that $\gamma\in \Gamma^*$ is non trivial. 
	Observe that, since the boundary of  $\D_\xi^{t}$ is a union of geodesic rays emanating from the point $x_\xi^{t}$, the geodesic segment $\mathopen]x_\xi^t, \gamma x_\xi^t\mathclose[$ stays out of both half-spaces $\D_\xi^{t}$ and $\gamma \D_\xi^{t}$.
\end{minipage}
\begin{minipage}{6cm}
	\hspace{0.1cm}	\includegraphics[width=5.4cm]{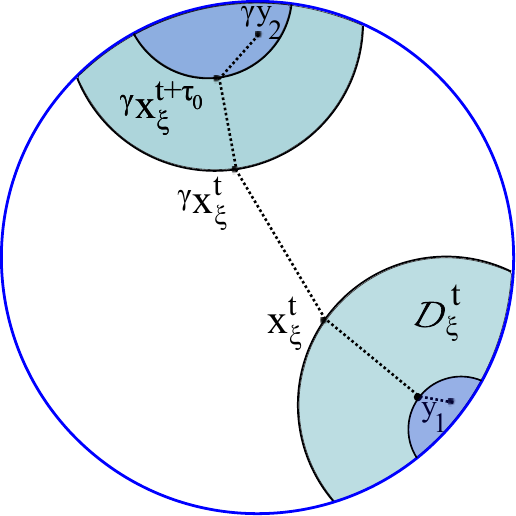}
\end{minipage}

\vspace*{0.8em}Let $\varepsilon=\inf\{ d(y,\gamma y)\; |\, y\in X,\; \gamma\in \Gamma^*\}>0$, and  $\alpha_\varepsilon$ be the corresponding angle in Lemma \ref{lem quad X}. 
As in the proof of Proposition \ref{prop a-s k-l}, choose $\tau_0$ large enough so that any half-space $\D_\xi^{t+\tau_0}$ is seen from the point $x_\xi^t$ under an angle less than $\a_\varepsilon$.  We may then apply Lemma \ref{lem quad X} to the quadrilateral $[y_1, x_\xi^t, \gamma x_\xi^t,\gamma y_2]\subset X$ to obtain
\begin{align*}
d(y_1,\gamma y_2)&\geq  d(y_1,x_\xi^t)+ d(x_\xi^t,\gamma x_\xi^t)   +d(\gamma x_\xi^t,\gamma y_2) -2\tilde\delta_{\varepsilon }\cr
&\geq  d(y_1,x_\xi^t) +d(y_2,x_\xi^t) -2\tilde\delta_{\varepsilon } \, .
\end{align*}	
Applying Lemma \ref{lem right angle triangle in X} to each triangle $[y_i,x_\xi^{t+\tau_0}, x_\xi^t]$ ($i=1,2$) yields
$$d(y_i,x_\xi^t)\geq d(y_i,x_\xi^{t+\tau_0})+\tau_0 -\hat\delta_0\, .$$
This is a contradiction if $\tau_0 > (\hat\delta_0 +\hat\delta_{\varepsilon })$, since  the triangle inequality yields
\begin{equation}
d(y_1,\gamma y_2)\geq d(y_1,y_2) +2\tau_0 -2(\hat\delta_0 +\hat\delta_{\varepsilon})\, . \plouf
\end{equation}

\subsection{Harmonic measures and embedded half-spaces in $ M$}\label{subsec:harm mes embedded half}
To construct an analogous to the Anderson-Schoen barrier functions for the quotient manifold $ M=\Gamma\backslash X$, we  work in the Hadamard manifold $X$. 

\vspace{0.3em}
Recall that we fixed a base point $o\in X$ and that  $\Omega_0\subset \Omega_\Gamma$ is relatively compact.
We keep the notation of Proposition \ref{prop a-s k-l} and Lemma \ref{lem t0 pour embedded}.
\bp\label{prop mesure anderson schoen} There exists $C_0\geq 1$  such that, for every $\xi\in \Omega_0$,
the harmonic measure at infinity of the saturation of ${\cal S}_\xi ^{t}$ under $\Gamma$ satisfies 
\begin{alignat*}{2}
\omega_o\bigl( \bigcup_{\gamma\in\Gamma}\gamma \, {\cal S}_\xi ^{t}\bigr)&\leq C_0\, e^{-\a t}& \quad&\text{for every $t\geq t_0$}
\cr\omega_{y}\bigl( \bigcup_{\gamma\in\Gamma}\gamma \, {\cal S}_\xi ^{t}\bigr)&\geq 1/C_0&\quad &\text{for every $y\in \D_\xi^{t+T}$.}
\end{alignat*}
\ep
We will need the following fact relative to the domain of discontinuity of $\Gamma$.
\bl\label{lem xi loin Lambda} Let $L\subset\Omega_\Gamma$ be any compact subset of the domain of discontinuity. Then, there is only a finite number of  elements $\gamma\in\Gamma^*$ with  $\gamma o\neq o$ and $\xi_{_-}\!(\gamma)\in L$.
\el
\proof Let $(\gamma_n)$ be a sequence of pairwise distinct elements of $\Gamma$ such that $\xi_{_-}\!(\gamma_n)\to\zeta\in\pX$. Since $\Gamma$ is discrete, $d(\gamma_n^{-1}o,o)\to\infty$ so that $\gamma_n^{-1}o\to\zeta$ hence $\zeta\in\Lambda_\Gamma$. The result follows readily.\eproof

\espace\noi{\bf Proof of Proposition \ref{prop mesure anderson schoen}}  The lower bound is an immediate consequence of Proposition \ref{prop a-s k-l}.

\vspace{0.4em}
Let now $\varepsilon>0$ small enough so that $\angle_o(\xi, \Lambda_\Gamma)\geq 2\varepsilon$ for every $\xi\in \Omega_0$. 
Lemma \ref{lem xi loin Lambda} ensures that the set $\Gamma_\varepsilon=\{ \gamma\in\Gamma^* \, |\; \angle_o(\xi_{_-}\!(\gamma), \Lambda_\Gamma)\geq\varepsilon\}$ is finite.

Assuming that $t\geq t_0$, we now seek an upper bound for $\omega_o( \bigcup_{\gamma\in\Gamma} \gamma \, {\cal S}_\xi ^{t})$.
We first observe that Corollary \ref{cor omega gammaA} and Lemma \ref{lem t0 pour embedded}
 ensure that
$$
e^{-c_\circ} \omega_o({\cal S}_\xi ^{t}) \;\bigl(\sum_{\gamma\in\Gamma} e^{-|\gamma|_o} \bigr) \leq \sum_{\gamma\in\Gamma}\omega_o( \gamma \, {\cal S}_\xi ^{t})=\omega_o( \bigcup_{\gamma\in\Gamma} \gamma \, {\cal S}_\xi ^{t})\, ,
$$
so that the series $\sum_{\gamma\in\Gamma} e^{-|\gamma|_o}$ converges.
When $\gamma\in\Gamma\smallsetminus \Gamma_\varepsilon$ is non trivial, one has $\angle_o(\xi, \xi_{_-}\!(\gamma))\geq\varepsilon$ for every $\xi\in \Omega_0$. Hence Corollary \ref{cor omega gammaA}  again yields
\begin{align*}
\sum_{\gamma\in\Gamma}\omega_o( \gamma \, {\cal S}_\xi ^{t}) &\leq \Bigl(e^{c_\circ} \, \sum_{\Gamma_\varepsilon}  e^{|\gamma|_o}+e^{a_\varepsilon} \, \sum_{\Gamma\smallsetminus \Gamma_\varepsilon} e^{-|\gamma|_o}\Bigr)\; \omega_o({\cal S}_\xi ^{t})\, ,
\end{align*}
and the claim now follows  from Proposition \ref{prop a-s k-l}.
\eproof

\espace\noi Recall that $V$ is the compact convex subset of $ M$ that we introduced in Proposition \ref{prop smooth core}. From now on, we will assume that the base point $o\in X$ is so chosen  that its projection $m_0\in M$ belongs to $V$.

\bc\label{cor upper bound harmonic} {\bf\small Upper bound for harmonic measures on $ M$}

\vsp There exists a constant $C_1$ such that the following holds. For every compact domain with smooth boundary $W$ whose interior contains $V$, every $\xi\in \Omega_0$, every $t\geq0 $ and every point $m\in V$~:
\begin{equation} \label{eq majo mesure}
\sigma_{m,W} (\partial W\cap D_\xi ^t) \leq C_1 e^{-\a t}\, .
\end{equation}
\ec
\proof 
It suffices to prove the assertion when $t$ is large. Recall that $T$ has been defined in Proposition \ref{prop a-s k-l}. For $\xi\in \Omega_0$ and $t\geq t_0+T$, introduce the positive harmonic function defined by
$$\tilde\eta _\xi^t : y\in X\to  \omega_{y}\bigl( \bigcup_{\gamma\in\Gamma}\gamma \, {\cal S}_\xi ^{t-T}\bigr)\in[0,1]\, .$$
The function $\tilde\eta _\xi^t $ is $\Gamma$-invariant and thus goes to the quotient to a harmonic function $\eta _\xi^t: M\to [0,1]$.  Proposition \ref{prop mesure anderson schoen} (that we apply to $t-T\geq t_0$) ensures that $\eta_\xi^t$ satisfies
\begin{alignat*}{2}
\eta _\xi^t(p)&\geq 1/C_0&\ &\text{for every $p\in D_\xi ^t$}\cr
\eta _\xi^t(m_0)&\leq C_0e^{\a T}\, e^{-\a t}\, .&&
\end{alignat*}
Applying the Harnack inequality (Lemma \ref{lem harnack}) to  $\tilde\eta _\xi^t $  yields
\begin{equation}\label{eq mino en m}
\eta _\xi^t(m)\leq e^{c_1d_V}C_0e^{\a T}\, e^{-\a t}\quad \hbox{for every $m\in V$,}
\end{equation}
where $d_V$ denotes the diameter of the compact convex set $V$.

\vspace{0.3em}
The function $C_0\, \eta _\xi^t$  is everywhere positive, and is greater or equal to $1$ on $\partial W\cap D_\xi^t$. Thus, the maximum principle ensures that 
$$
C_0\, \eta _\xi^t (p)\geq \sigma_{p,W}( \partial W\cap D_\xi ^t)
$$
holds for every $p\in W$, and the claim now follows from \eqref{eq mino en m}.
\eproof

\espace
Note that the  constants $t_0$, $\tau_0$, $C_0$ and $C_1$ that we introduced in the previous paragraphs  depend only on the group $\Gamma$, on the base point $o$, on  $\Omega_0$ and on  the compact subset $V\subset M$.

\subsection{Gromov products and embedded half-spaces}\label{subsec:gromov parallel}
We end this chapter with Proposition \ref{prop deux demi espaces}, that relates half-spaces corresponding to the same point at infinity  with level sets of Gromov products.

\vspace{0.6em}
Recall that we choose a base point $o\in X$ whose projection $m_0\in M$ lies in the compact subset $V\subset M$ and
that  $t_0$ and $\tau_0\geq 1$ were defined in Lemma \ref{lem t0 pour embedded} and Proposition \ref{prop tau0}.

\bp\label{prop deux demi espaces}  
There exists a constant $g$ such that,
for every $\xi\in \Omega_0$, every $t\geq t_0+\tau_0$, 
every point $p\in D_\xi ^{t+\tau_0 }$ and every point $m\in V$~:
\begin{align*}
\{   q\in M \ | \    (p,q)_{m}\geq t+ g \} \subset D_\xi ^t\, .
\end{align*}
\ep

\proof
Let $p\in D_\xi^{t+\tau_0 }$ and $q\notin D_\xi ^t$. 
Consider a minimizing geodesic segment $[p,q]\subset  M$ and introduce two points $p'$ and $q'$ where $[p,q]$ intersects the hyperplanes $H_\xi^ {t+\tau _0}$ and $H_\xi ^t$. 

Since $t\geq t_0+\tau_0$, it follows from Proposition \ref{prop tau0} that any, hence the only, minimizing quadrilateral $[p',m_\xi^{t+\tau_0},m_\xi^t,q']$ lies in the embedded half-plane $D_\xi^{t_0}$ and is thus isometric to a quadrilateral in $X$.  

\noindent	\begin{minipage}{6.4cm}
	This quadrilateral $[p',m_\xi^{t+\tau_0},m_\xi^t,q']$ is right-angled at both vertices $m_\xi^t$ and $m_\xi^{t+\tau_0 }$, and $d(m_\xi^t,m_\xi^{t+\tau_0 })=\tau_0\geq 1$. 	Hence Lemma \ref{lem quad X}  applies to prove  that the point $ m_\xi^t$ is within distance $\tilde\delta_1$ of  the edge $[p', q']$.  Since  $[p',q']\subset [p,q]$ and $d(m,m_\xi^t)\leq t+d_V$,  it follows from Lemma \ref{lem proprietes Gromov prdct}  that 
	$$(p,q)_{m}\leq d(m,[p,q])\leq d(m,[p',q'])\leq
	t+(d_V+\tilde\delta_1)\, .$$
\end{minipage}
\begin{minipage}{7.5cm}
	\hspace{0.3em}	\includegraphics[width=6cm]{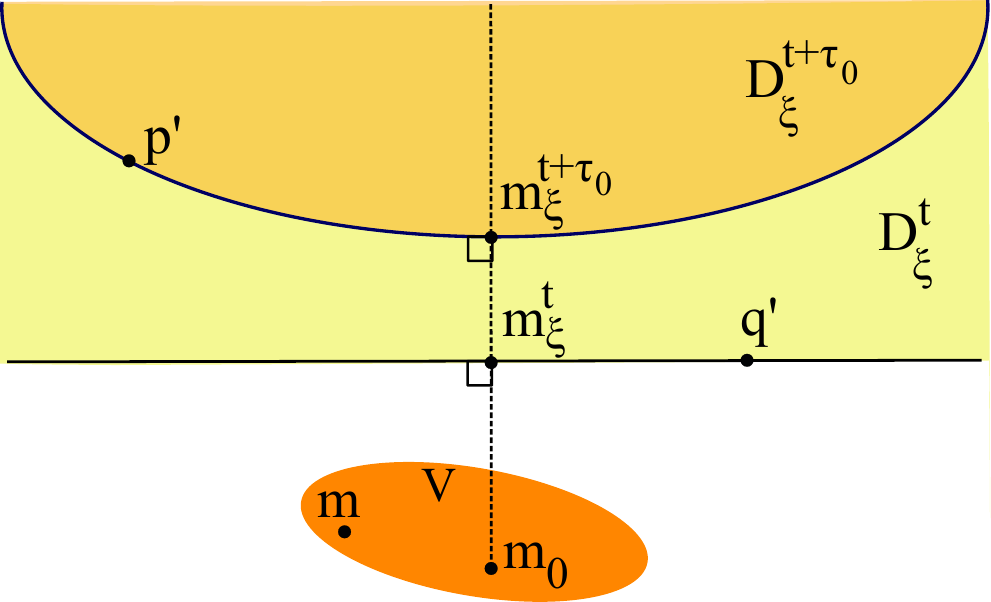}
\end{minipage}

\vspace{0.2em}\noi
The claim follows for $g >d_V+\tilde\delta_1$.
\eproof

\section{Interior estimates}\label{sec:interior estimates}

\begin{quote}
	In this final section, we wrap up the proof of Proposition \ref{prop rho}, that gives a uniform bound for the distances $d(f,h_R)$. 
\end{quote}
	We split the proof  into two parts. 
In the first part, where we assume that the point $m\in V_R$ where the distance $d(f,h_R)$ is reached  lies far away from the convex core, the proof reduces  to the proof of the main theorem of \cite{BH18}.

In the second part, where we assume that the point $m$ lies in a fixed neighbourhood of the convex core, we must deal with the topology of the quotient $ M=\Gamma\backslash X$.
\subsection{Harmonic quasi-isometric maps $H:X\to Y$}\label{subsec: XY}
In \cite{BH18}, we proved that a quasi-isometric map $F:X\to Y$ between two pinched Hadamard manifolds $X$ and $Y$ is within bounded distance from a unique harmonic map. As in the present paper, this harmonic map was obtained as the limit of a family of solutions of Dirichlet problems on bounded domains with boundary value $F$, where a uniform bound for the distances between $F$ and the solutions of these Dirichlet problems ensured the convergence of the family. 

In the following technical statement, which is local in nature, we gather some information obtained in \cite{BH18} that was used to obtain this uniform bound. The first part of our proof of Proposition \ref{prop rho} will derive easily from this statement, see  Proposition \ref{prop far}.
\bfa \label{fact}   Let $c\geq 1$. There exist  $\ell_0>1$ and $\rho_0>0$ with the following property.

Let  $F, H:B(x,\ell_0)\to Y$ be two smooth maps defined on a ball $B(x,\ell_0)\subset X$ with radius $\ell_0$, and such that the distance $$\displaystyle \rho:=\sup_{z\in B(x,\ell_0)} d(F(z),H(z))$$  is reached at the center $x$ of the ball, namely  $\rho =d(F(x),H(x))$. Assume that $H$ is harmonic and that the map $F$ satisfies
\begin{align}
\label{eq F qi}	\text{$c^{-1} d(z,z')-c\leq d(F(z),F(z'))\leq c\, d(z,z')$ for any $z,z'\in B(x,\ell_0)$,} 
\end{align}
Then
$$\rho\leq\rho_0\, .$$
\efa
\proof This fact is proven in \cite[Section 4]{BH18}. Alternatively, one may  follow  paragraphs \ref{subsec: overview close} through \ref{subsec:estim gromov 0} below, replacing the domain $V_{\ell}$ with the ball $B(m,\ell_0)$  and using the uniform estimates \eqref{eq estimees uniformes mesure X} for the harmonic measures of balls in the Hadamard manifold $X$ instead of the new estimates \eqref{eq minoration mesure} and \eqref{eq majo mesure}. \eproof

\subsection{Part one~: estimate far away from the core}\label{subsec:far}
In this paragraph, we introduce a finite family of embedded half-spaces in $ M$, whose union is a neighbourhood of infinity, and that will be used throughout the proof of Proposition \ref{prop rho}. Then, we prove Proposition \ref{prop rho} in case the distance $d(f,h_R)$ is reached far away from the convex core.

\vspace{0.8em}
Recall that, after smoothing, the map $f: M\to Y$ is assumed to be $c$-Lipschitz \eqref{eq f smooth} and that each point $\xi\in\partial_\infty M$ admits a neighbourhood to which $f$ restricts as a quasi-isometric map. 

From now on,  $\Omega_0\subset\Omega_\Gamma$ is a fixed compact neighbourhood of a fundamental domain for the action of $\Gamma$ on $\Omega_\Gamma$.

\bl\label{lem choix de V} (1) There exists a finite number of embedded half-spaces  $D_{\zeta_i}^{t_i}$ ($1\leq i\leq N$) with $\zeta_i\in \Omega_0$ and $t_i\geq t_0$ such that, taking perhaps a larger constant $c$  in \eqref{eq f smooth}~:
\begin{itemize}
	\item each restriction $f:D_{\zeta_i}^{t_i}\to Y$ is a quasi-isometric map with constant $c$
	\item $\cup_{i=1}^N D_{\zeta_i}^{t_i+\tau_0}$  is a neighbourhood of infinity in $M$.
\end{itemize}
(2) We may also assume  that $V$ has been choosen large enough so that
\begin{align}
\label{eq V dans les demi}& \bullet \;\partial V\subset \cup_{i=1}^N D_{\zeta_i}^{t_i+\tau_0}\\
\label{eq rayon inj}& \bullet \; \hbox{for $m\notin V$, one has $\inj (m)\geq \ell_0$}\\
\label{eq qi boules} & \bullet \; \hbox{for $m\notin V$, the restriction of $f$ to  $B(m,\ell_0)$ is $c$ quasi-isometric.}
\end{align}
\el
\proof  (1)  We proved in Lemma \ref{lem t0 pour embedded} that, for every $\xi\in\Omega_0$ and any $t\geq t_0$, the half-space $D_\xi^t$ embeds in $ M$. Hence, the claim follows from the hypothesis on $f$, since the boundary at infinity $\partial_\infty  M$ is compact.

(2) The injectivity radius $\inj : M\to\plu$ is a proper function, and the  complement in $ M$ of
$$
\cup_{i=1}^N\{ m\in M \; |\; B(m,\ell_0)\subset D_{\zeta_i}^{t_i}\}
$$
is bounded. 
Hence, it suffices to replace the convex set $V$ of Proposition \ref{prop smooth core} by its $R_0$-neighbourhood  for some large  $R_0$ to ensure  these conditions.
\eproof 

\espace
We want a uniform upper bound for the distance $d(f,h_R)$ when $R$ is large. We may thus assume that 
\begin{equation}\label{eq R ell0}R\geq \ell_0+1\, .\end{equation}
In the next proposition, we obtain such a bound in case the distance $d(f,h_R)$ is reached at some point $m$ which is far away from the convex core, that is if $m\notin V$. The  case where $m\in V$ will be carried out in the next paragraphs.

\bp\label{prop far} Suppose that the distance  $d(f,h_R)=d(f(m),h_R(m))$
is reached at some point $m\in V_R\smallsetminus V$. Then
$$
d(f,h_R)\leq\rho_0+B\ell_0\, ,
$$
where $B$ is the constant in Proposition \ref{prop estimee dist bord}, and $\rho_0$ is defined in Fact \ref{fact}.
\ep
\proof Assume first that $d(m,\partial V_R)\leq \ell_0$.   It follows from Proposition \ref{prop estimee dist bord} that  
$$d(f,h_R)=d(f(m),h_R(m))\leq B\, d(m,\partial V_R)\leq B\ell_0 \, .$$
Assume now that $d(m,\partial V_R)> \ell_0$, so that
the ball $B(m,\ell_0)$ lies in $V_R$. Since $m\notin V$,  Condition \eqref{eq rayon inj} ensures that this ball $B(m,\ell_0)$ is isometric to a ball with radius $\ell_0$ in the Hadamard manifold $X$. Fact \ref{fact} applies to the restrictions $F=f_{|B(m,\ell_0)}$ and $H=(h_R)_{|B(m,\ell_0)}$, so that $\rho\leq\rho_0$.
The result follows. \eproof

\subsection{Part two~: estimate close to the core, an overview}\label{subsec: overview close}

To complete the proof of Proposition \ref{prop rho}, we  assume from now on that the distance $\rho=d(f,h_R)$ is reached at some point $m$ that belongs to the fixed compact set $V$. We pick a large $\ell$ (namely $\ell$ will have to satisfy Conditions \eqref{eq ell1 delta}, \eqref{eq ell1 ell(s)} and \eqref{eq ell1 trop grand final}), and we will mainly work on  the compact convex set with smooth boundary $V_{\ell}$ which is the $\ell$-neighbourhood of $V$. Note that $V_\ell$  does not depend on $R$. 

\espace
We will assume that
\begin{equation}\label{eq ell1 delta}
\ell\geq d_V +1\, ,
\end{equation}
where $d_V$ denotes the diameter of the domain $V$.
Since we want an upper bound for  the distance $\rho =d(f,h_R)$ when $R$ is large, we may also assume from now on that 
\begin{align}\label{eq rho c}
\rho&\geq c \\
\label{eq R ell1} R&\geq \ell+2\, .
\end{align}
For any point $p$ in the $1$-neighbourhood of $V_\ell$, namely for $p\in V_{\ell+1}$, Condition \eqref{eq R ell1} ensures that $B(p,1)\subset V_R$, so that Corollary \ref{cor cheng} yields
\begin{align}\label{eq cheng interior estimates}
\| D_ph_R\|\leq\kappa\rho\, .
\end{align}
We introduce $y=f(m)\in Y$, which is the image under $f$ of the point $m\in V$ where the distance $d(f,h_R)$ is reached.  For any point $p\in\partial V_\ell$, we shall
study the three following Gromov products relative to this point $y$~:  
\begin{equation*}
g_0(p)=(f(p),h_R(m))_y , \ g_1(p)=(f(p),h_R(p))_y, \ g_2(p)=(h_R(p),h_R(m))_y\, .
\end{equation*}
\begin{center}
\includegraphics[height=5cm]{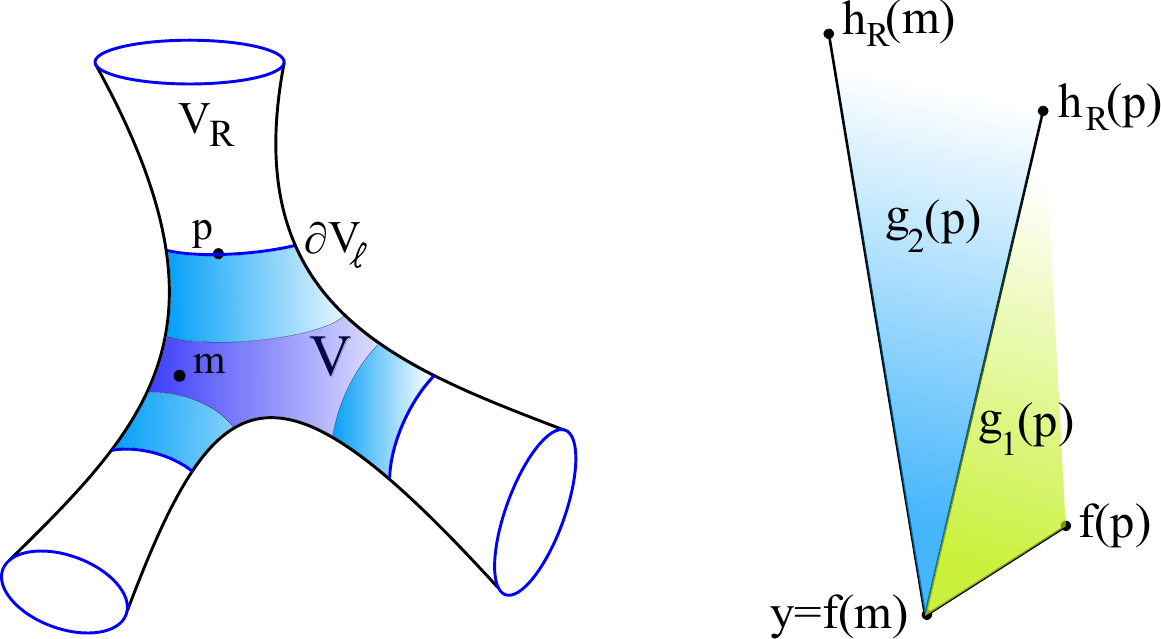}
\end{center}
If $\rho$ is large, we shall prove that on a suitable subset $U_{\ell, R}$ of the boundary $\partial V_\ell$,  both $g_1$ and $g_2$ are large (Lemma \ref{lem g1 sur U} and Corollary \ref{cor g2 sur U}) 
while the measure of $U_{\ell,R}$ is large enough (Lemma \ref{lem measure U}) to ensure that $g_0$ cannot be that large on the whole $U_{\ell,R}$ (Lemma \ref{lem g0 sur U}). This will yield a  contradiction thanks to Inequality \eqref{eq gromov product inegalite triangulaire} satisfied by Gromov products.

The arguments we develop here are similar to those of \cite[Section 4]{BH18}, and that led to Fact \ref{fact}. In the setting of our previous paper \cite{BH18}, they relied on the uniform upper and lower bounds for  harmonic measures of balls in the Hadamard manifold relative to their center. In our new context, they rely on the uniform upper and lower bounds for  harmonic measures obtained in Proposition \ref{prop lower bound harmonic} and Corollary \ref{cor upper bound harmonic}.

\subsection{The subset $U_{\ell,R}\subset\partial V_\ell$}\label{subsec:U}
	We introduce the domain $U_{\ell,R}\subset\partial V_\ell$ that will play a central role in the proof of Proposition \ref{prop rho}.
\bd\label{def U}
Let $U_{\ell,R}$ be the set of those points  $p\in\partial V_\ell$ where the distance $d(y,h_R(p))$ is close to $\rho=d(y,h_R(m))$, namely
\begin{equation*}
U_{\ell,R}=\{ p\in \partial V_\ell\, |\; d(y,h_R(p))\geq\rho-\frac{\ell}{2c}\}\, .
\end{equation*}
\ed
In the next lemma, we give a lower  bound for the ``size" of the domain $U_{\ell,R}\subset \partial V_\ell$, which is uniform with respect to $R$ and to the choice of $\ell$. 
\bl\label{lem measure U}
The harmonic measure $\sigma_{m,V_\ell}$ of the set $U_{\ell,R}\subset \partial V_\ell$ relative to the point $m\in V$ satisfies
$$
\sigma_{m,V_\ell} (U_{\ell,R})\geq\, \frac 1{5\, c^2}.
$$
\el
\proof  
We first observe that, for any point $q\in V_\ell$,  one has
\begin{equation}
d(y,h_R(q))\leq \rho + 2c\, \ell\, .\label{eq u bdd}
\end{equation}
Indeed, since by  \eqref{eq f smooth} the function $f$ is $c$-Lipschitz, the triangle inequality yields
$$
d(f(m),h_R(q))\leq d(f(m),f(q))+d(f(q),h_R(q)) \leq c\, (d_V+\ell)+\rho
$$
where as usual $d_V$ denotes the diameter of $V$. Thus \eqref{eq u bdd} follows, since we assumed in  \eqref{eq ell1 delta} that $\ell\geq d_V$. 

\vspace{0.3em}
Now, Lemma \ref{lem distance subharmonic} asserts that the function $u:q\to d(y,h_R(q))-\rho$ is subharmonic on $V_R$. Moreover  \eqref{eq u bdd} ensures that $u$ is bounded above by $ 2c\, \ell$  on $V_\ell$. Since $u(m)=0$, this yields
$$
0 \leq\int_{\partial V_\ell} u(p)\, d\sigma_{m,V_\ell} (p)\leq \sigma_{m,V_\ell}(U_{\ell,R})\, ( 2c\, \ell)-(1-\sigma_{m,V_\ell} (U_{\ell,R}))\, \frac{\ell}{2c}\, ,
$$
thus $1\leq \sigma_{m,V_\ell} (U_{\ell,R}) \, (1+4c^2)$ hence the claim, since we assumed  $c\geq 1$.\eproof

\subsection{Upper bound for the Gromov product $g_0(p)$ on $U_{\ell,R}$}\label{subsec:estim gromov 0}
	In this paragraph we prove  that, if $\ell$ is large enough, the Gromov products  $g_0(p)=(f(p),h_R(m))_y$ cannot be uniformly large	on the whole $U_{\ell,R}$. 

\vspace{0.6em}
We first prove that the image $f(U_{\ell,R})\subset Y$, seen from the point $y=f(m)$, is relatively spread out. 
\bl\label{lem g0 sur U}  There exists a distance $\bar\ell$ and a constant $\bar g_0>0$ such that, if 
\begin{equation}
\label{eq ell1 ell(s)}
\ell\geq \bar\ell \,  ,
\end{equation}
then there exist  two points $p_1,p_2\in U_{\ell,R}\subset \partial V_\ell$ with
\begin{equation*}
(f(p_1),f(p_2))_{y}\leq \bar g_0\, .
\end{equation*}
\el
\bc\label{cor majoration g0} If $\ell\geq\bar\ell$, then there exists a point $p\in U_{\ell,R}$ such that
$$
g_0(p)\leq \bar g_0+2\delta_Y\, .
$$
\ec
\proof 
Follows from Lemma \ref{lem g0 sur U} and Inequality \eqref{eq gromov product inegalite triangulaire}.
\eproof

\espace\noi{\bf Proof of Lemma \ref{lem g0 sur U}} We first construct the points $p_1,p_2\in U_{\ell,R}$.

We proved in Lemma \ref{lem measure U} that the harmonic measure $\sigma_{m,V_\ell}$ of $U_{\ell,R}$ is bounded below by $1/(5c^2)$.
It thus follows from \eqref{eq V dans les demi} that, for any choice of $\ell$, there exists an index $1\leq j\leq N$ with $\sigma_{m,V_\ell}(U_{\ell,R}\cap D_{\zeta_j}^{t_j+\tau_0})\geq 1/(5Nc^2)$.

Fix $t$ such that  $C_1e^{-\a t}<1/(5Nc^2)$. Assume moreover that $t\geq t_0+\tau_0$. Let $\bar \ell >0$ be large enough so that 
$\partial V_{\ell}\subset  \cup_{\xi\in \Omega_0} D_\xi ^{t+\tau_0}$ for each $\ell\geq\bar\ell$.
Pick a point $p_1\in U_{\ell,R}\cap D_{\zeta_j}^{t_j+\tau_0}$ and choose $\xi\in \Omega_0$ such that $p_1\in D_\xi^{t+\tau_0 }$. 
Proposition \ref{prop deux demi espaces}, that applies since  $t\geq t_0+\tau_0$,  and Corollary \ref{cor upper bound harmonic} ensure that
\begin{align*} 
\sigma_{m,V_\ell} (\{ p\in U_{\ell,R}\, |\; (p_1,p)_{m}\geq t+ g\})\leq \sigma_{m,V_\ell} (D_\xi^t) \leq C_1 e^{-\a t}<1/(5Nc^2)\, .
\end{align*}
Hence, there exists $p_2\in U_{\ell,R}\cap D_{\zeta_j}^{t_j+\tau_0}$ such that $(p_1,p_2)_{m}\leq t+ g$. 

\vspace{0.8em} We now turn our attention to the two images $f(p_1), f(p_2)$.
If  $f: M\to Y$ were  assumed to be  quasi-isometric with constant $c$,  we would infer immediately from Lemma \ref{lem gromov product quasi iso}  that
$(f(p_1),f(p_2))_{y}\leq \bar g_0$, with $\bar g_0=c\, (t+g)+A$.

	But under the hypotheses of Theorem \ref{th main}, where we only assume that the restriction of $f$ to the half-space $D_{\zeta_j}^{t_j}$ is a $c$ quasi-isometric map, we have to make a slight adjustment to this elementary proof. For simplicity of notation, let us denote by $m_j=m_{\zeta_j}^{t_j+\tau_0}$.
	Introduce the compact set $W=V\cup\{m_1,\cdots , m_N\}$, and let $d_W$ be its diameter. Observe that
$$
(p_1,p_2)_{m_j} \leq (p_1,p_2)_{m}  +d_W\leq t+g+d_W\, .
$$

We are now ready to use the fact that the restriction of $f$ to $D_{\zeta_j}^{t_j}$ is a quasi-isometric map. Indeed, the three points $p_1,p_2,m_j$ belong to the half-space $D_{\zeta_j}^{t_j+\tau_0}$, whose convex hull is included in $D_{\zeta_j}^{t_j}$ (Proposition \ref{prop tau0}). This convex hull being isometric to a convex subset of $X$, Lemma \ref{lem gromov product quasi iso}  yields
$$
(f(p_1),f(p_2))_{f(m_j)} \leq c\, (t+ g+d_W)+A\, .
$$
To prove our claim,  replace the origin $f(m_j)$ with $y=f(m)$ in this Gromov product, observing that, since $f$ is $c$-Lipschitz, we have
$d(f(m_j),f(m))\leq c\, d_W$.
\eproof

\subsection{Lower bound for the Gromov product $g_1(p)$ on $U_{\ell,R}\subset\partial V_\ell$}\label{subsec:estim gromov 1}
Our second estimate for the Gromov products is the only one that relies on the left-hand side of Condition \eqref{def qi}.

\bl\label{lem g1 sur U} 
There exists a constant $\bar g_1$ such that
\begin{equation*}
g_1(p)=(f(p),h_R(p))_y\geq\frac{\ell}{4c}-\bar g_1\, 
\end{equation*}
holds for every $p\in U_{\ell,R}$.
\el
\proof In case the map $f$ is $c$ quasi-isometric on the whole $ M$, the proof is straightforward. Indeed, using the bound $d(f(p),h_R(p))\leq\rho$, the definition of $U_{\ell,R}$  and observing that, since $m\in V$ and $p\in\partial V_{\ell}$, one has $d(m,p)\geq\ell$, we obtain
\begin{align*}
2(f(p),h_R(p))_y &=d(f(p),f(m))+d(h_R(p),f(m))-d(f(p),h_R(p)) \cr
&\geq (\frac{\ell}c-c) +(\rho-\frac{\ell}{2c}) -\rho = \frac{\ell}{2c}-c\, .
\end{align*}

When the map $f$ is not supposed to be globally  quasi-isometric, we proceed  as in the proof of Lemma \ref{lem g0 sur U} and introduce the index  $1\leq j\leq N$ such that $p\in D_j$, so that
$$
d(f(p),f(m_j))\geq (1/c)\, d(p,m_j)-c \, .
$$
The same computation as above gives
$$
2(f(p),h_R(p))_ {f(m_j)} \geq \frac\ell{2c} -c-2c\, d_W\, .
$$ Using again the fact that $f$ is $c$-Lipschitz on $M$ to change the base point, we obtain the result.
\eproof

\subsection{Lower bound for the Gromov product $g_2(p)$ on $\partial V_\ell$}\label{subsec:estim gromov 2}

This last estimate relies on the uniform lower bound for the harmonic measures of a family of suitable subdomains of $ M$ proven in   Paragraph \ref{subsec: harmonic measure quotient}.
\bl\label{lem rho/2}
There exists a constant $\rho_2(\ell)>0$ that depends only on $\ell$ (and not on $R$) and such that if
\begin{equation}\label{eq rho rho1}
\rho>\rho_2(\ell)
\end{equation}
then
\begin{equation*}
d(y,h_R(q))\geq \rho /2
\end{equation*}
holds for every point $q\in V_\ell$.
\el
\proof Assume that there exists a point $q\in V_\ell$ where  $d(y,h_R(q))< \rho /2$. Because of the bound \eqref{eq cheng interior estimates} for the covariant derivative of $h_R$ on $V_{\ell+1}$, it follows that one has  $d(y,h_R(z))\leq 3\rho /4$ for every point $z$ in the  ball $B(q,1/4\kappa )$.

\vspace{0.6em} 
Let $\lambda =1/4\kappa$ and $L=2\ell$, so that \eqref{eq ell1 delta} yields $0< d(m,q)\leq L$. 
With the notation of Proposition \ref{prop lower bound harmonic},  we introduce the constant $s_2(\ell)=s(\lambda , L)$.

We proceed as in the proof of Lemma \ref{lem measure U}.  Consider the subharmonic function $u:z\to d(y,h_R(z))-\rho$ on the domain  $W_{m,q}$  we introduced in Proposition \ref{prop lower bound harmonic}. This function $u$ vanishes at the point $m$. 
We just proved that $u\leq -\rho /4$ on the ball $B(q,\lambda)$, while thanks to \eqref{eq ell1 delta}~:
$$
u(z)\leq d(f(m),f(z))+d(f(z),h_R(z))-\rho\leq 2c\, \ell
$$
for any point $z\in V_{\ell+1}$, and in particular on the boundary $\partial W_{m,q}$. We thus infer that
$0 \leq 	-\rho \, \sigma_{m,W_{m,q}} (B(q,\lambda))+8c\,\ell$,
hence $\rho\leq {8c\ell}/{s_2(\ell)}$. This proves our claim, with $\rho_2(\ell)=8c\ell/{s_2}(\ell)$.
\eproof

\bc\label{cor g2 sur U} 
There exists a constant $\bar g_2(\ell)$, that depends on $\ell$ but not on $R$, such that  if $\rho>\rho_2(\ell)$
\begin{equation*}
g_2(p)=(h_R(m),h_R(p))_y\geq\frac\rho 2-\bar g_2(\ell)\, \log\rho \, 
\end{equation*}
holds for any point $p\in \partial V_\ell$.
\ec
\proof Let $p\in\partial V_\ell$ and let $[m,p]$ be a minimizing geodesic segment from $m$ to $p$.  By Assumption \eqref{eq ell1 delta}, its length is at most $d_V+\ell\leq 2\ell$.
We infer from the bound \eqref{eq cheng interior estimates} for the covariant derivative of $h_R$ that the length of the curve $h_R([m,p])\subset Y$ is at most $2\ell\kappa\rho$. 

Since this curve stays away from the large ball $B(y,\rho /2)$, it will look short seen from the point $y$. Indeed,  select a subdivision $(z_i)_{0\leq i\leq 2^n-1}$ of  $h_R([m,p])$, with  $ 1\leq d(z_i,z_{i+1})\leq 2$. One thus has $n\leq {\rm log}_2\rho +\nu$ for some constant $\nu >0$.
Since $d(z_{2i},z_{2i+1})\leq 2$, Lemma \ref{lem rho/2} gives $(z_{2i},z_{2i+1})_y\geq\rho /2-1$ when $0\leq i\leq 2^{n-1}$.
Then, the triangle inequality for Gromov products (Lemma \ref{lem proprietes Gromov prdct}) yields
$(z_{2i},z_{2(i+1)})_y\geq \rho /2 -1- \delta_Y$ when $0\leq i\leq 2^{n-1}$. Iterating the process yields $(h_R(m),h_R(p))_y \geq\frac\rho 2 -1 -n\, \delta_Y$ as claimed. \eproof

\subsection{Proof of Proposition \ref{prop rho}} \label{subsec:proof majo rho}
	We now prove that, if $\rho$ is too large, the estimates for the three Gromov products $g_0$, $g_1$ and $g_2$ that we obtained in the previous sections lead to a contradiction, thus completing the proof of Proposition \ref{prop rho}  and of our main theorem \ref{th main}.

\vspace{0.6em}
Let us first stress the fact that  both constants $\bar g_0$ and $\bar g_1$ do not depend on $R$, $\rho$ nor $\ell$, while $\bar g_2$ depends on $\ell$ but not  on $R$ nor $\rho$.

We begin by choosing a radius $\ell$ large enough to satisfy   \eqref{eq ell1 delta} and  \eqref{eq ell1 ell(s)}, as well as
\begin{equation}\label{eq ell1 trop grand final}
\frac{\ell}{4c}-\bar g_1-4\delta_Y >\bar g_0\, ,
\end{equation}
where $\bar g_0$ and $\bar g_1$ are defined in Lemmas  \ref{lem g0 sur U} and \ref{lem g1 sur U}.

We let $R\geq \ell_0+\ell+2$, so that Conditions \eqref{eq R ell0} and \eqref{eq R ell1} are satisfied.  Assume by contradiction that the  distance $\rho:=\rho_R=d(f,h_R)$ is very large, namely that $\rho$ satisfies \eqref{eq rho c}, \eqref{eq rho rho1} and 
\begin{equation}\label{eq rho plus que ell1 gromov}
\frac\rho 2-\bar g_2(\ell)\log\rho\geq \frac{\ell}{4c}-\bar g_1\, .
\end{equation}
We proved in Lemma \ref{lem g1 sur U} and Corollary \ref{cor g2 sur U} that, under these assumptions~:
\begin{alignat*} {3}
(f(p),h_R(p))_y&=g_1(p)\geq\frac{\ell}{4c}-\bar g_1& \qquad  &\text{when $p\in U_{\ell,R}$} \\
(h_R(m),h_R(p))_y&=g_2(p)\geq\frac\rho 2-\bar g_2(\ell)\, \log\rho &\qquad& \text{when $p\in\partial V_\ell$.}
\end{alignat*}
Thus  \eqref{eq gromov product inegalite triangulaire} and \eqref{eq rho plus que ell1 gromov} yield, for any point $p\in U_{\ell,R}$, the lower bound
\begin{align*}
(f(p),h_R(m))_y &=g_0(p)\geq\min(g_1(p),g_2(p))-2\delta_Y\cr
&\geq \frac{\ell}{4c}-\bar g_1-2\delta_Y> \bar g_0+2\delta_Y
\end{align*}
thanks to our choice of $\ell$ in  \eqref{eq ell1 trop grand final}. This is a contradiction to  Corollary \ref{cor majoration g0}. This ends the proof of Proposition \ref{prop rho} in case $m\notin V$. 

\vspace{0.3em}\noi
The case where $m\in V$ has already  been dealt with in Proposition \ref{prop far}.
\eproof

\vfill\eject
\small{
\bibliographystyle{plain}
\bibliography{harmonic3}
}
\vspace{2em}

{\small\noindent Y. Benoist  \& D. Hulin,\; 
	CNRS \& Université Paris-Saclay \\ Laboratoire de mathématiques d’Orsay, 91405, Orsay, France
	
	\vspace{0.3em}
\noindent yves.benoist@u-psud.fr \;\&\; dominique.hulin@u-psud.fr}
\end{document}